\documentclass[11pt,a4paper,twoside]{article}
\usepackage{geometry}
\usepackage[english]{babel}
\usepackage[utf8]{inputenc}
\usepackage[small,it,textfont={small,it}]{caption}
\usepackage{amsmath,amssymb,amsfonts,amsthm}

    \DeclareMathOperator{\tr}{tr}
    \newcommand{\exn}[2]{{#1}\mathrm{e}{#2}}
\usepackage{cite} % Declare this if using Bibtex
\usepackage{graphicx}
\usepackage{enumitem}
    \newlist{abbrv}{itemize}{1}
    \setlist[abbrv,1]{label=,labelwidth=1cm,align=parleft,
    itemsep=0.1\baselineskip,leftmargin=!}
\usepackage{hyperref}
\usepackage{url}
\usepackage{textcomp}
\usepackage{xspace}
\usepackage[toc,page]{appendix}
\usepackage[dvipsnames]{xcolor}
\usepackage{array}
\usepackage{booktabs}
\newcommand{\bfg}[1]{\boldsymbol{#1}} %Bold style for symbols
\newcommand{\MATLAB}{\textsc{Matlab}\xspace}
\title{\textbf{Stability Assessment of Stochastic Differential-Algebraic Systems via Lyapunov
Exponents with an Application to Power Systems}}
\author{Andr\'es Gonz\'alez-Zumba\footnotemark[1] \,
and Pedro Fern\'andez-de-C\'ordoba\footnotemark[2] \,
and\\ Juan-Carlos Cort\'es\footnotemark[3] \,
and Volker Mehrmann\footnotemark[4]}
\date{\today}
\begin{document}
\setlength{\parindent}{12pt}
\setlength{\parskip}{8pt}
\maketitle
%%%%%%%%%%%%%%%%%%%%%%%%% ABSTRACT %%%%%%%%%%%%%%%%%%%%%%%%%%%%%%%%%%%%%%%%%%%%%%%%%%%%%%%%%%%%%%%%
\begin{abstract}
\noindent
In this paper we discuss Stochastic Differential-Algebraic Equations (SDAEs) and the asymptotic
stability assessment for such systems via Lyapunov exponents (LEs). We focus on  index-one SDAEs
and their reformulation as ordinary  stochastic differential equation (SDE). Via ergodic theory it
is then feasible to analyze the LEs via the random dynamical system generated by the underlying
SDE. Once the existence of well-defined LEs is guaranteed, we proceed to the use of numerical
simulation techniques to determine the LEs numerically. Discrete and continuous $QR$
decomposition-based numerical methods are implemented to compute the fundamental solution matrix
and to use it in the computation of the LEs. Important computational features of both methods are
illustrated via numerical tests. Finally, the methods are applied to two applications from power
systems engineering, including the single-machine infinite-bus (SMIB) power system model.
\end{abstract}
%%%%%%%%%%%%%%%%%%%%%%%%%%%%%%%%%%%%%%%%%%%%%%%%%%%%%%%%%%%%%%%%%%%%%%%%%%%%%%%%%%%%%%%%%%%%%%%%%%%

%\noindent
%{\bf AMS subject classification.:} XXXXX, XXXXX, XXXXX.

\maketitle
\renewcommand{\thefootnote}{\fnsymbol{footnote}}
\footnotetext[1]{Departamento de Matem\'atica Aplicada, Universitat Polit\`ecnica
de Val\`encia, Camino de Vera s/n 46022, Valencia, Spain. \texttt{jorgonzu@posgrado.upv.es}.}
\footnotetext[2]{Instituto Universitario de Matem\'atica Pura y Aplicada, Universitat Polit\`ecnica
de Val\`encia, Camino de Vera s/n 46022, Valencia, Spain. \texttt{pfernandez@mat.upv.es}.}
\footnotetext[3]{Instituto de Matem\'atica Multidisciplinar, Universitat Polit\`ecnica
de Val\`encia, Camino de Vera s/n 46022, Valencia, Spain. \texttt{jccortes@mat.upv.es}.}
\footnotetext[4]{Institut f\"ur Mathematik MA 4-5, Technische Universit\"at Berlin, Str. des 17.
Juni 136, D-10623
Berlin, FRG. \texttt{mehrmann@math.tu-berlin.de}.}
\renewcommand{\thefootnote}{\arabic{footnote}}

%%%%%%%%%%%%%%%%%%%%%%%%%%%%%%%%%%%%%%%%%%%%%%%%%%%%%%%%%%%
%%%%%%% SECTION 1: INTRODUCTION %%%%%%%%%%%%%%%%%%%%%%%%%%%
%%%%%%%%%%%%%%%%%%%%%%%%%%%%%%%%%%%%%%%%%%%%%%%%%%%%%%%%%%%
\section{Introduction}\label{sec:intro}
Modeling the dynamic behavior of systems employing differential-algebraic equations
%coupled to algebraic equations, which represent the constraints of such systems,
is a mathematical representation paradigm widely used in many areas of science and engineering.
%These equations are well known as differential-algebraic equations (DAEs).
On the other hand, the dynamics of systems perturbed by stochastic processes have been adequately
modeled by stochastic differential equations (SDEs). The need of a generalized concept which covers
both DAEs and SDEs, and allows the modeling and analysis of constrained systems subjected to
stochastic disturbances, has led to the formulation of stochastic differential-algebraic equations
(SDAEs). While there has been a broad and fruitful development both in the field of DAEs and SDEs
(see, e.g. \cite{BieCM12,BreCP96,KunM06} and \cite{KloP10,Mao07,Oks13}, respectively), studies of
the concepts and the numerical treatment of SDAEs are rather limited, see e.g.
\cite{ConT10,SchD98,Win03}. The main reason for this is that the proper treatment of algebraic
constraints in SDAEs faces many difficulties, except in the case that all the constraints are
explicitly given and can be resolved during the numerical integration process, which is the case
that we will discuss.

As tool for the stability analysis we use Lyapunov exponents (LEs) introduced as characteristic
exponents in \cite{Lya92}. The theory of LEs experienced a crucial development with the work
\cite{Ose68} which, via the Multiplicative Ergodic Theorem (MET), ensures the regularity and
existence of the LEs belonging to a linear cocycle over a metric dynamical system, see  \cite{Arn98}
for a detailed  presentation of the theory. We review and extend  the main concepts of this
approach for asymptotic stability assessment of differential-algebraic equations driven by Gaussian
white noise and apply the technique in the setting of power systems.

%We start by surveying the theory of SDAEs. In particular, we are study autonomous
%SDAE systems of the form
%
%\[
%\bfg Ed\bfg x_t = \bfg f_0(\bfg x_t)dt + \sum_{j=1}^{m} \bfg f_j(\bfg x_t) dw_t^j, \qquad
%t\in\mathbb{R}^+,
%\]
%
%together with a consistent initial value $\bfg x_{t_0}=\bfg x_0$.
Following the ideas from \cite{ConT10,ConT12,KupKR12,SchD98,Win03}, properties such as the
existence and uniqueness of solutions are reviewed. Analogously to the DAE case, we define the
class of strangeness-free (index-one) SDAEs and show that SDAE systems with this structure can be
reduced to a classical SDE system that describes the dynamics of the original SDAE for which the
MET can be applied
%Additionally, we know from the classical
%theory that an autonomous It\^o SDE generates a RDS (see \cite{Arn98}). Therefore, the Oseledets'
%MET can be applied
to define the LEs of the system.

Once we have extended the theoretical framework that guaratees the existence of well-defined LEs,
we study numerical methods, based on the $QR$ factorization of the fundamental solution matrix,
that allow the numerical computation of spectral values asociated to the Lyapunov spectra. The
first technique requires computing the fundamental solution matrix and forming an orthogonal
factorization, while the second one involves performing a continuous $QR$ decomposition of the
fundamental solution matrix. Both techniques have been extensively studied in deterministic ODE and
DAE systems, see \cite{BenGGS80a,BenGGS80b,DieV03,DieV07,LinM09,LinMV11}. This paper follows the
ideas exposed in \cite{CarBJ10}, where these $QR$ methods were extended to the stochastic case.

Finally, these concepts and computational techniques are used to assess the asymptotic stability of
power systems affected by stochastic fluctuations. We illustrate the results with elementary test
cases such as a single-machine infinite-bus system.
% This is an interesting attemp to evidence
%the usefulness of the LEs as an indicator of stability for dynamical systems working urder
%uncertainty and open new avenues to its potential application in the stability analysis of
%large-scale systems.

The paper is organized as follows: In Section \ref{sec:theory} we recall the main theory of
strangeness-free SDAEs, their relation with the SDEs, and the existence of LEs % of the RDS
generated by such SDEs. Section \ref{sec:lenum} presents the discrete and continuous $QR$-based
decomposition methods and their evaluation. % through a single example.
Interesting study-cases in the power systems area are presented in Section \ref{sec:apps}. We
finish with some conclusions in Section \ref{sec:conclu}.

%%%%%%%%%%%%%%%%%%%%%%%%%%%%%%%%%%%%%%%%%%%%%%%%%%%%%%%%%%%
%%%%%%% SECTION 2: THEORETICAL REVIEW %%%%%%%%%%%%%%%%%%%%%
%%%%%%%%%%%%%%%%%%%%%%%%%%%%%%%%%%%%%%%%%%%%%%%%%%%%%%%%%%%
\section{Review of the theory}\label{sec:theory}
\subsection{Stochastic Differential-Algebraic Equations}\label{sec:sdae}
%differential equations on manifolds.
Consider a system of quasi-linear stochastic differential-algebraic equations  (SDAEs) of the form
\begin{equation} \label{eq:ito-nlSDAE}
    \bfg Ed\bfg x_t = \bfg f_0(\bfg x_t)dt + \sum_{j=1}^{m} \bfg f_j(\bfg x_t) dw_t^j, \quad
    t\in\mathbb{I}:=[t_0,t_f],
\end{equation}
with a singular matrix $\bfg E\in\mathbb{R}^{n\times n}$ of rank $d<n$.
%is known as \textit{``leading matrix''}.
%has the form $E = \bigl[\begin{smallmatrix} E_{11}&0\\ 0&0 \end{smallmatrix}\bigr]$,
The function $\bfg f_0\in\mathcal{C}^k(\mathbb{D}_{\bfg x},\mathbb{R}^n)$ (for some $k\geq1$) is
known as \emph{drift}, and $\bfg f_1,\ldots,\bfg f_m\in\mathcal{C}^{k+1}(\mathbb{D}_{\bfg
x},\mathbb{R}^n)$ are the \emph{diffusions}, here $\mathbb{D}_{\bfg x}\subseteq\mathbb{R}^n$
is an open set. Furthermore, $w_t^j$ (for $j=1,\ldots,m$) form an $m$-dimensional \emph{Wiener
process} defined on the complete probability space $(\Omega,\mathcal{F},\mathbb{P})$ with a
filtration $(\mathcal{F}_t)_{t\geq t_0}$. Each $j$-th Wiener process is understood as a process with
independent increments such that $(w_t-w_s)\sim\mathcal{N}(0,t-s)$, i.e. is a Gaussian random
variable for all $0\leq s<t$, such that
\[
w_0=0, \quad \mathbb{E}[w_t-w_s]=0, \quad \mathbb{E}[w_t-w_s]^2=t-s.
\]

Since the Wiener process is nowhere differentiable, it is more convenient to represent equation
(\ref{eq:ito-nlSDAE}) in its integral form as
\begin{equation} \label{eq:ito-nlSDAE-if}
	\bfg E\bfg x_t = \bfg E\bfg x_{t_0} + \int_{t_0}^t \bfg f_0(\bfg x_s)ds + \sum_{j=1}^{m}
    \int_{t_0}^t \bfg f_j(\bfg x_s) dw_s^j, \quad t\in\mathbb{I}.
\end{equation}
Here, the first integral is a stochastic Riemann-Stieltjes integral, and the second one is a
stochastic It\^o-type integral, see e.g. \cite{Oks13}.
%In this regard, two approaches of stochastic integrals are known:
%the It\^o's and the Stratonovich's approaches. In this paper, we will adopt the It\^o's
%interpretation. However, swapping between It\^o and Stratonovich interpretations is feasible by
%converting the SDAEs through a simple formula that relates the two definitions.

%We assume initial values $\bfg x_{t_0}=\bfg x_0\in L^2(\Omega)$
%$\bfg x_t\in\mathcal{L}^2(\Omega,\mathcal{F},\mathbb{P};\mathbb{R}^n)$ are unknown variables.
We assume consistent initial values $\bfg x_{t_0}=\bfg x_0$ independent of the Wiener processes
$w_t^j$ and with finite second moments \cite{Oks13}. A solution $\bfg x_t = \bfg x(t,\omega)$ of
(\ref{eq:ito-nlSDAE-if}) is an $n$-dimensional vector-valued Markovian stochastic process
depending on $t\in\mathbb{I}$ and $\omega\in\Omega$ (the parameter $\omega$ is commonly omitted in
the notation of $\bfg x$). Such a solution can be defined as \emph{strong solution} if it fulfills
the following conditions, see e.g. \cite{ConT12,Win03}.
\begin{itemize}
    \item $\bfg x(\cdot)$ is adapted to the filtration $(\mathcal{F}_t)_{t\geq t_0}$,
    %(i.e., it does not depend on future information),
    \item $\int_{t_0}^{t_f}|f_0^\ell (x_s)|ds<\infty$ almost sure (a.s.), for all
    $\ell=1,\ldots,n$,
    \item $\int_{t_0}^{t_f}|f_j^\ell (x_s)|^2dw_s^j<\infty$ a.s., for all $ j=1,\ldots,m$, and
    $\ell=1,\ldots,n$,
    \item (\ref{eq:ito-nlSDAE-if}) holds for every $t\in\mathbb{I}$ a.s.
    %and the integrals exist and the the equations are satisfied in (\ref{eq:ito-nlSDAE-if}) for
    %every $t\in\mathbb{I}$ a.s.
\end{itemize}

Because of the presence of the algebraic equations associated with the kernel of $\bfg E$, the
solution components associated with these equations would be directly affected by white noise and
not integrated. In order to avoid this, a reasonable restriction is to ensure that the noise sources
do not appear in the algebraic constraints. According to \cite{SchD98,Win03}, this assumption can
be accomplished in SDAE systems whose deterministic part
\begin{equation}\label{eq:nlDAE}
	\bfg E\dot{\bfg x}_t = \bfg f_0(\bfg x_t), \quad t\in\mathbb{I},
\end{equation}
is a DAE with tractability index one \cite{LamMT13,Win03}, in which the constraints are regularly
and globally uniquely solvable for parts of the solution vector.
%Based on these ideas, we state the following definition:
%
%\begin{definition}(\cite{Win03})\label{def:SDAE_i-1}
%An SDAE system of the form (\ref{eq:nlDAE}) is called an index-one SDAE if
%the noise sources do not appear in the algebraic constraints, and the constraints can be solved
%globally and uniquely for a subset of the variables.
%\end{definition}
We slightly modify this assumption and consider SDAE systems whose deterministic part
(\ref{eq:nlDAE}) is a regular \textit{strangeness-free} DAE \cite{KunM06} i.e. it has
differentiation index one. A system with these characteristics can be transformed into a
semi-explicit form by means of an appropriate kinematic equivalence transformation
\cite{LinM09,BieCM12}, i.e., there exist pointwise orthogonal matrix functions
$\bfg{\mathcal{P}}$ and $\bfg{\mathcal{Q}}$ such that, pre-multiplying (\ref{eq:ito-nlSDAE}) by
$\bfg{\mathcal{P}}$, and changing the variables $\bfg x_t$ according to the transformation $\bfg
x_t=\bfg{\mathcal{Q}}\hat{\bfg x}_t$
%\bfg{\mathcal{P}}\bfg E\bfg{\mathcal{Q}}= [I 0,0 0]
one obtains a system in semi-explicit form
\begin{subequations}\label{eq:nlSDAEsp}
    \begin{align}
        d \hat{\bfg x}_t^D &= \hat{\bfg f}_0^D(\hat{\bfg x}_t^D,\hat{\bfg x}_t^A)dt + \sum_{j=1}^{m}
        \hat{\bfg f}_j^D(\hat{\bfg x}_t^D,\hat{\bfg x}_t^A)dw_t^j,\label{eq:nlSDEsp}\\
        0 &= \hat{\bfg f}_0^A(\hat{\bfg x}_t^D,\hat{\bfg x}_t^A)dt + \sum_{j=1}^{m} \hat{\bfg
        f}_j^A(\hat{\bfg x}_t^D,\hat{\bfg x}_t^A)dw_t^j,\label{eq:nlAEsp}
    \end{align}
\end{subequations}
where $\hat{\bfg x}_t^D$ and $\hat{\bfg x}_t^A$ is a separation of the transformed state into
differential and algebraic variables, respectively, that is performed in such a way that the
Jacobian of the function $\hat{\bfg f}_{0}^A$ with respect to the algebraic variables is
nonsingular, see \cite{KunM06} for details of the construction. The condition that the noise sources
do not appear in the constraints, implies that $\sum_{j=1}^{m} \hat{\bfg f}_j^A\equiv0$, so that
the algebraic equations in (\ref{eq:nlAEsp}) can be solved as $\hat{\bfg x}_t^A=\bfg F^A(\hat{\bfg
x}_t^D)$ and inserted in the dynamic equations (\ref{eq:nlSDEsp}) yielding an ordinary SDE
\begin{align}\label{eq:uSDE}
    d \hat{\bfg x}_t^D &= \hat{\bfg f}_0^D(\hat{\bfg x}_t^D,\bfg F^A(\hat{\bfg x}_t^D))dt +
    \sum_{j=1}^{m}\hat{\bfg f}_j^D(\hat{\bfg x}_t^D,
    \bfg F^A(\hat{\bfg x}_t^D))dw_t^j.
%    \quad     \textrm{with} \qquad \hat{\bfg x}_t^A = \bfg F^A(\hat{\bfg x}_t^D).
\end{align}
This equation is called \emph{underlying SDE} of the strangeness-free SDAE. It acts in the
lower-dimensional subspace $\mathbb{R}^d$, with $d=n-a$ (where $a$ denotes the number of algebraic
equations). The SDE system (\ref{eq:uSDE}) preserves the inherent dynamics of a strangeness-free
SDAE system~\cite{LamMT13}. Note that in this way, the algebraic equations have been removed from
the system, but whenever a numerical method is used for the numerical integration, then one has to
make sure that the algebraic equations are properly solved at each time step, so that the
back-transformation to the original state variables can be performed.

\subsection{Random Dynamical Systems generated by SDEs}\label{sec:rds}
In the previous section we have discussed the reduction of an autonomous strangeness-free SDAE to
its underlying SDE, which preserves the dynamic characteristics of the original system. Using the
back-transformation the definitions and properties attributed to the underlying SDE, and the
analysis performed on it can be extended to the original  SDAE. For simplicity we use the following
representation, where  the drift and diffusion terms  are combined into one term.
\begin{equation} \label{eq:ito-SDE}
    d\bfg x_t=\bfg f_0(\bfg x_t)dt+\sum_{j=1}^{m} \bfg f_j(\bfg x_t)dw_t^j =
    \sum_{j=0}^{m}\bfg f_j(\bfg x_t)dw_t^j,\quad t\in\mathbb{I},
\end{equation}
where
$\bfg f_0\in\mathcal{C}_b^{k,\delta}$, $\bfg f_1,\dots,\bfg f_m\in\mathcal{C}_b^{k+1,\delta}$ and
$\sum_{j=1}^{m}\sum_{i=1}^{d} f_j^i\frac{\partial}{\partial x_i}f_j\in\mathcal{C}_b^{k,\delta}$
for some $k\geq1$ and $\delta>0$. Here $\mathcal{C}_b^{k,\delta}$ is the Banach space of
$\mathcal{C}^k$ vector fields on $\mathbb{R}^d$ with linear growth and bounded derivatives up to
order $k$ and the $k$-th derivative is $\delta$-H\"{o}lder continuous. In addition, we assume that
the differential operator $L:=\bfg f_0+\frac{1}{2}\sum_{j=1}^m (\bfg f_j)^2$ is strong hypoelliptic
in the sense that the Lie algebra $\mathcal{L}(\bfg f_0, \bfg f_1,\ldots, \bfg f_m)$ generated by
the vector fields $\bfg f_j$ (with $j=0,\ldots,m$) has dimension $d$ for all $\bfg
x_t\in\mathbb{R}^d$~\cite{Arn98}. Once again, $w_t^j$ (for $j=0,\ldots,m$) is a $m$-dimensional
Wiener process, this time with the convention $dw_t^0\equiv dt$.

For a given initial value, the solution process generates a \emph{Markovian stochastic process},
and the SDE (\ref{eq:ito-SDE}) generates a \emph{random dynamical system (RDS)}
$\varTheta=(\theta,\bfg\varphi)$ which is an object consisting of a \emph{metric dynamical system
(MDS)} $\theta$ for modeling the random perturbations, and a \emph{cocycle}
$\bfg\varphi:\mathbb{R}^+\times\Omega\times\mathbb{R}^d\rightarrow\mathbb{R}^d$ over this
system. The ergodic MDS is denoted by
$\theta\equiv(\Omega,\mathcal{F},\mathbb{P},(\theta_t)_{t\in\mathbb{R}})$ with the filtration
$(\mathcal{F}_t)_{t\geq t_0}$, and defined by the \emph{Wiener shift}
\[
\theta_t\omega(\cdot)= \omega(t+\cdot)-\omega(t), \quad t\in\mathbb{I},
\]
which means that a shift transformation given by $\theta$ is measure-preserving and ergodic
\cite{CarH17}.

%Thm 2.3.26 (Pag.83), Thm 2.3.32 (Pag.93), and Thm 2.3.39 (Pag.98)
Together with the SDE (\ref{eq:ito-SDE}), we define the variational system
\begin{equation}\label{eq:ito_veSDE}
    d\bfg v_t = \sum_{j=0}^{m} \bfg{J}_j(\bfg x_t)\bfg v_t \ dw_t^j, \quad t\in\mathbb{I},
    \quad \textrm{with} \quad \bfg{J}_j(\bfg x_t):=\left(\frac{\partial \bfg f_j}{\partial
    \bfg x}\right)\bigg\rvert_{\bfg x(t;t_0)},
\end{equation}
obtained after linearizing (\ref{eq:ito-SDE}) along a solution. If we denote by
$\bfg{\varPhi}(t,\omega,\bfg x)$ the Jacobian of $\bfg\varphi(t,\omega)$ at $\bfg x_t$, then
$\bfg{\varPhi}$ is the unique solution of the variational equation (\ref{eq:ito_veSDE}), satisfying
\begin{equation}\label{eq:sol_veSDE}
    \bfg{\varPhi}(t,\omega,\bfg x) = \bfg I_d + \sum_{j=0}^m\int_0^t \bfg{J}_j
    (\bfg\varphi(s)\bfg x)\bfg{\varPhi}(s,\bfg x)dw_s^j,\quad t\in\mathbb{I},
\end{equation}
where $\bfg I_d$ denotes the identity matrix of size $d$. Therefore, $\bfg{\varPhi}$ it is a matrix
cocycle over $\varTheta$. The system (\ref{eq:ito-SDE})-(\ref{eq:ito_veSDE}) uniquely generates a
$\mathcal{C}^{k-1}$ RDS $(\bfg\varphi,\bfg{\varPhi})$ over $\theta$. Moreover, the determinant of
$\bfg{\varPhi}$ satisfies the \emph{Liouville equation}
\begin{equation}\label{eq:liouville}
    \det \bfg{\varPhi} = \exp \left(\sum_{j=0}^m\int_0^t \tr
    (\bfg{J_j}(\bfg\varphi(s)\bfg x)dw_s^j \right),\quad t\in\mathbb{I},
\end{equation}
being then a scalar cocycle over $\varTheta$ (see Theorems 2.3.32 and 2.3.39-40 in \cite{Arn98} for
details).

\subsection{Lyapunov Exponents of Ergodic RDSs}\label{sec:met}
An important result of the theory of RDSs is the so-called Multiplicative Ergodic Theorem (MET)
developed in \cite{Ose68}. This concept allows the definition of LEs for linear
cocycles over a ergodic MDS. First, the MET assumes for the linear cocycle $\bfg{\varPhi}$ that the
integrability condition
\[
\log^+\lVert \bfg{\varPhi}(t,\omega,\bfg x)\rVert\in\mathcal{L}^1,
\]
is satisfied, where $ \log^+$ denotes the positive part of $\log$. This guarantees that the
variational equation (\ref{eq:ito_veSDE}) associated with (\ref{eq:ito-SDE}) is well-posed.
Additionally, let $\mu$ be an ergodic invariant measure with respect to the cocycle
$\bfg\varphi$~\cite{Arn98}. Then, the MET assures the existence of an invariant set
$\hat\Omega\subset\Omega$ of full $\mu$-measure, such that for each $\omega\in\hat\Omega$ there
is a measurable decomposition
\[
\mathbb{R}^d = L_1(\omega)\oplus\dots\oplus L_p(\omega),
\]
of $\mathbb{R}^d$ into random linear subspaces $L_{i}(\omega)$, which are invariant under
$\varTheta$. Here $p\leq d$, where $d_i\in\mathbb{N}$ denotes the dimension of the subspace
$L_i(\omega)$ (with $1\leq i\leq p$), and $\sum_{i=1}^p d_i=d$. This splitting is dynamically
characterized by real numbers $\lambda_1>\ldots>\lambda_p$ which quantify the exponential growth
rate of the subspaces. These are called \emph{Lyapunov exponents}, and are defined by
\[
\lambda_i:=\lim_{t\to\infty}\frac{1}{t}\log\lVert\bfg{\varPhi}(t,\omega,\bfg
x)\rVert\quad \mbox{\rm whenever} \quad \bfg x\in L_i(\omega)\setminus\{0\}.
\]
According to \cite[pag.~118]{Arn88}, the LEs $\lambda_i$ are independent of $(\omega,\bfg x)$ and
thus they are universal constants of the cocycle generated by (\ref{eq:ito_veSDE}) under the
ergodic invariant probability measure $\mu$. Finally, the following identity holds for
$\bfg{\varPhi}$ when the system is Lyapunov \emph{regular}~\cite{PikP16,CarBJ10}
\begin{equation}\label{eq:reg}
    \sum_{i=1}^d \lambda_i = \lim_{t\to\infty}\frac{1}{t}\log|\det \bfg{\varPhi}(t,\omega,\bfg x)|.
\end{equation}
In practice, it is hard (if not impossible) to verify Lyapunov regularity for a particular system
\cite{Arn98}. One of the key statements of the MET is that linear RDS (whether these are constant,
periodic, quasi-periodic, or almost-periodic) are a.s. Lyapunov regular.

The concept of LEs plays an important role in the asymptotic stability assessment of dynamical
systems subjected to stochastic disturbances. Under appropriate regularity assumptions, the
negativity of all LEs of the system of variational equations implies the exponential asymptotic
stability of both the linear SDE and the original nonlinear SDE system.

%%%%%%%%%%%%%%%%%%%%%%%%%%%%%%%%%%%%%%%%%%%%%%%%%%%%%%%%%%%
%%%%%%% SECTION 3: METHODS FOR COMPUTING LES %%%%%%%%%%%%%%
%%%%%%%%%%%%%%%%%%%%%%%%%%%%%%%%%%%%%%%%%%%%%%%%%%%%%%%%%%%
\section{$\bfg{QR}$ Methods for computing LEs}\label{sec:lenum}
In this section we derive the numerical techniques to compute the finite-time approximation of the
LEs. In the same way as \cite{CarBJ10}, this paper proposes an adaptation of the ideas from
\cite{DieRV97,DieV03,DieV07} to the stochastic case (linear RDSs). The methods take advantage on
the existence of a Lyapunov transformation of the linear RDS to an upper-triangular structure, and
the feasibility to retrieve a numerical approximation of the LEs from that form. The transformation
is performed through an orthogonal change of variables. The approach is made under the assumption of
Lyapunov regularity of the system. In order to explain the methods, let us consider again the
SDE as an initial value problem of the form
\begin{equation}\label{eq:SDE3}
    d\bfg x_t = \sum_{j=0}^{m} \bfg f_j(\bfg x_t)dw_t^j, \quad t\in\mathbb{I}, \quad
    \bfg x_{t_0}=\bfg x_{t_0},
\end{equation}
where $\bfg f_j$ are sufficiently smooth functions. The corresponding  variational equation of
(\ref{eq:SDE3}) along with the solutions $\bfg x_t(\bfg x_0)$, turned into a matrix initial value
problem, is given by
\begin{equation}\label{eq:veSDE}
    d\bfg V_t = \sum_{j=0}^{m} \bfg J_j(\bfg x_t) \bfg V_t dw_t^j, \qquad \bfg V_0=\bfg I_d,
\end{equation}
with the identity matrix $\bfg I_d\in\mathbb{R}^{d\times d}$ as initial value,  where
$\bfg{J}_j(\bfg x_t):=\frac{\partial \bfg f_j}{\partial\bfg x}$ are the Jacobians of the vector
functions $\bfg f_j(\bfg x_t)$, and $\bfg V\in\mathcal{C}^1(\mathbb{I}\times\mathbb{R}^{d\times
d})$ is the fundamental solution matrix, whose columns are linearly independent solutions of the
variational equation. A key theoretical tool for determining the LEs is the computation of the
continuous $QR$ factorization of $\bfg V_t$,
\[
\bfg V_t= \bfg Q_t\bfg R_t,
\]
where $\bfg Q_t$ is orthogonal, i.e., $\bfg Q_t^T\bfg Q_t=\bfg I_d$, and $\bfg R_t$ is upper
triangular with positive diagonal elements $R_t^{ii}$ for $i=1,\ldots,d$. Applying the MET theory
presented in Subsection~\ref{sec:met}, and taking into account the norm-preserving property of
the orthogonal matrix function $\bfg Q_t$, we have
\begin{equation} \label{eq:cle}
    \lambda_i = \lim_{t\to\infty}\frac{1}{t}\log\lVert \bfg V_t\bfg p_i\rVert =
    \lim_{t\to\infty}\frac{1}{t}\log \lVert \bfg R_t \bfg p_i\rVert,
\end{equation}
where ${\{\bfg p_i}\}$ is an orthonormal basis associated with the splitting of $\mathbb{R}^d$.
Lyapunov regular systems preserve their regularity under kinematic similarity transformations.
Then, considering the regularity condition (\ref{eq:reg}), the Liouville equation
(\ref{eq:liouville}), and performing some algebraic manipulations (see details in
\cite[pag.~150]{CarBJ10}), the LEs are given by
\begin{equation}\label{eq:ler}
    \lambda_i = \lim_{t\to\infty} \frac{1}{t} \log |R_t^{ii}| \quad \textrm{a.s.}, \quad
    \textrm{for} \quad i=1,\ldots,d.
\end{equation}
The $QR$ methods require to perform the $QR$ decomposition of $\bfg V_t$ for a long enough time, so
that the $R_t^{ii}$ have started to converge. Depending on whether the decomposition is performed
after or before integrating numerically the variational equation, the method is called
\emph{discrete or continuous $QR$ method}.

\subsection{Discrete $\bfg{QR}$ Method}\label{sec:cm}
The discrete $\bfg{QR}$ method is a very popular method for computing LEs in ODEs and DAEs. In this
approach, the fundamental solution matrix $\bfg V_t$ and its triangular factor $\bfg R_t$ are
indirectly computed by a reorthogonalized integration of the variational equation (\ref{eq:veSDE})
through an appropriate $QR$ decomposition. Thus, given grid points
$0=t_0<t_1<\ldots<t_{N-1}<t_N=T$, we can write $\bfg V_{t_\ell}$ in terms of the state-transition
matrices as
\begin{equation}\label{eq:dqr1}
    \bfg V_{t_\ell}=\bfg Z_{(t_\ell,t_{\ell-1})}\bfg Z_{(t_{\ell-1},t_{\ell-2})}
    \cdots \bfg Z_{(t_2,t_1)}\bfg Z_{(t_1,t_0)}\bfg V_{t_0}.
\end{equation}
At $t_0=0$, we perform a standard matrix $QR$ decomposition
\[
\bfg V_{t_0}=\bfg Q_{t_0}\bfg R_{t_0},
\]
%
%This implies $\bfg Q_{t_0}=\bfg I_d$ and $\bfg R_{t_0}=\bfg I_d$.
and for $\ell=1,2,\dots,N$, we determine $\bfg Z_{(t_\ell,t_{\ell-1})}$ as the numerical solution
(via numerical integration) of the matrix initial value problem
\begin{equation}\label{eq:D_QR_VE1}
    d \bfg Z_{(t_\ell,t_{\ell-1})}=\sum_{j=0}^{m}\bfg J_j(\bfg x_t)
    \bfg Z_{(t_\ell,t_{\ell-1})} dw_t^j, \quad
    \bfg Z_{(t_{\ell},t_{\ell-1})} = \bfg Q_{t_{\ell-1}}, \quad t_{\ell-1}\le t \le t_\ell,
\end{equation}
and then compute the $QR$ decomposition
\[
\bfg Z_{(t_{\ell},t_{\ell-1})}=\bfg Q_{t_{\ell}}\bfg R_{(t_\ell,t_{\ell-1})},
\]
where $\bfg R_{(t_\ell,t_{\ell-1})}$ has positive diagonal elements. From (\ref{eq:dqr1}), the
value of the fundamental matrix $\bfg V_{t_\ell}$ is determined via
\[
\bfg V_{t_\ell}=\bfg Q_{t_{\ell}}\bfg R_{(t_\ell,t_{\ell-1})}\bfg R_{(t_{\ell-1},t_{\ell-2})}
\cdots \bfg R_{(t_2,t_1)}\bfg R_{(t_1,t_0)}\bfg R_{t_0},
\]
which is again a $QR$ factorization with positive diagonal elements. Since this is unique, for the
$QR$ decomposition $\bfg V_{t_\ell}=\bfg Q_{t_\ell}\bfg R_{t_\ell}$, we have
\[
\bfg R_{t_\ell} = \bfg R_{(t_\ell,t_{\ell-1})}\bfg R_{(t_{\ell-1},t_{\ell-2})}
\cdots \bfg R_{(t_2,t_1)}\bfg R_{(t_1,t_0)}\bfg R_{t_0} = \prod_{\kappa=0}^\ell \bfg R_\kappa.
\]
Here we denote as $\bfg R_\kappa$ the triangular transition matrices
$\bfg R_{(t_\ell,t_{\ell-1})}$ with $\kappa=0,1,\ldots,\ell$. From (\ref{eq:ler}), the LEs are thus
computed  as
\begin{equation}\label{eq:D_QR_LEs}
	\lambda_i=\lim_{\ell\to\infty}\frac{1}{t_\ell}\log
	\left|\prod_{\kappa=0}^\ell R_\kappa^{ii}\right|
	= \lim_{\ell\to\infty}\frac{1}{t_\ell}\sum_{\kappa=0}^\ell\log|R_\kappa^{ii}|, \quad
    i=1,\ldots,d.
\end{equation}
\subsection{Continuous $\bfg{QR}$ Method}\label{sec:cm}
The implementation of the continuous $QR$ technique requires to determine a system of SDEs for the
$\bfg Q$ factor and the scalar equations for the logarithms of the diagonal elements of the $\bfg
R$ factor elementwise. Then, once the orthogonal matrix $\bfg Q$ is computed by numerical
integration, the logarithms of the diagonal elements of $\bfg R$ can also be obtained.

By differentiating in the It\^o sense the decomposition $\bfg V_t=\bfg Q_t \bfg R_t$ and using the
orthogonality $\bfg Q_t^T\bfg Q_t=\bfg I_d$, we get
\begin{align}
    d\bfg V_t &= (d\bfg Q_t)\bfg R_t+\bfg Q_t(d\bfg R_t),\label{eq:qrd}\\
    \bfg 0 &= (d\bfg Q_t^T)\bfg Q_t+\bfg Q_t^T(d\bfg Q_t).\label{eq:qqd}
\end{align}
Inserting (\ref{eq:qrd}) into the variational equation (\ref{eq:veSDE}), and multiplying by $\bfg
Q_t^T$ from the left and by $\bfg R_t^{-1}$ from the right, we obtain
\begin{equation}\label{eq:cqr1}
    \bfg Q_t^T(d\bfg Q_t) + (d\bfg R_t)\bfg R_t^{-1} = \sum_{j=0}^m \bfg Q_t^T \bfg J_j(\bfg x_t)
    \bfg Q_t dw_t^j.
\end{equation}
Since $(d\bfg R_t)\bfg R_t^{-1}$ is upper triangular, the skew-symmetric matrix $d\bfg S_t := \bfg
Q_t^T(d\bfg Q_t)$ satisfies
\begin{equation}\label{eq:cqr2}
    d\bfg S_t^{il} =
    \begin{cases}
    \ \sum_{j=0}^m \left(\bfg Q_t^T\bfg J_j(\bfg x_t)\bfg Q_t\right)^{jl}dw_t^j, & \quad i>l,\\
    \ 0, & \quad i=l,\\
    \ -\sum_{j=0}^m \left(\bfg Q_t^T\bfg J_j(\bfg x_t)\bfg Q_t\right)^{jl}dw_t^j, & \quad i< l.
    \end{cases}
\end{equation}
This results in an SDE for $\bfg Q_t$ given by
\begin{equation}\label{eq:cqr3}
    d\bfg Q_t = \bfg Q_t d\bfg S_t =\sum_{j=0}^m \bfg Q_t \bfg T_t^j(\bfg x_t,\bfg Q_t)dw_t^j,
\end{equation}
where the matrices $\bfg T_t^j(\bfg x_t,\bfg Q_t)$ (for $j=0,\ldots,m$) are defined via
\begin{equation}
    \left(\bfg T_t^j(\bfg x_t,\bfg Q_t)\right)^{il} =
    \begin{cases}
    \ \left(\bfg Q_t^T\bfg J_j(\bfg x_t)\bfg Q_t\right)^{jl}, & \quad i> l,\\
    \ 0, & \quad i=l,\\
    \ -\left(\bfg Q_t^T\bfg J_j(\bfg x_t)\bfg Q_t\right)^{jl}, & \quad i<l.
    \end{cases}
\end{equation}
A corresponding SDE for $\bfg R_t$ can be obtained from (\ref{eq:cqr1}) and (\ref{eq:cqr2}) via
\begin{equation}
    d\bfg R_t = \sum_{j=0}^m (\bfg Q_t^T \bfg J_j(\bfg x_t)\bfg Q_t-\bfg T_t^i(\bfg x_t,\bfg
    Q_t))\bfg R_tdw_t^j,
\end{equation}
and the equation for the $i$th diagonal element $R_t^{ii}$ is given by
\begin{equation}\label{eq:cqr4}
    d R_t^{ii} = \sum_{j=0}^m (\bfg Q_t^T \bfg J_j(\bfg x_t)\bfg Q_t)^{ii}
    R_t^{ii}dw_t^j, \quad \textrm{for } i=1,\ldots,d.
\end{equation}
Since the computed LEs can be obtained from (\ref{eq:ler}), we make use of the It\^o Lemma to
introduce the following SDE for the function $\psi_t^i=\log R_t^{ii}$ from (\ref{eq:cqr4}),
\begin{equation}\label{eq:cqr5}
	d\psi_t^i = d(\log R_t^{ii}) = \sum_{j=0}^m(\bfg Q_t^T \bfg J_j(\bfg x_t)\bfg
    Q_t)^{ii}dw_t^j -\frac{1}{2}\left[\sum_{j=0}^m(\bfg Q_t^T \bfg J_j(\bfg x_t)\bfg
    Q_t)^{ii}dw_t^j \right]^2.
\end{equation}
If we assume that there are no correlations between the diffusion terms in the SDE system, then we
do not have terms $dw_t^k\, dw_t^\ell$ (for $1<k<m$, and $1<\ell<m$, with $k\neq \ell$) in the
SDE (\ref{eq:cqr5}).% , when the square polynomial is developed.
Also, using that  $dt\, dt\equiv 0$,
$dt\, dw_t^k\equiv 0$, and $dw_t^j\, dw_t^j\equiv dt$ for $1<k<m$, the SDE (\ref{eq:cqr5})
is reduced to
\begin{equation}\label{eq:cqr6}
	d\psi_t^i = \sum_{j=0}^m(\bfg Q_t^T \bfg J_j(\bfg x_t)\bfg
    Q_t)^{ii}dw_t^j -\frac{1}{2}\sum_{j=1}^m\left[(\bfg Q_t^T \bfg J_j(\bfg x_t)\bfg
    Q_t)^{ii}\right]^2dt.
\end{equation}
By integrating this SDE, it is possible to obtain the LEs $\lambda_i$ from
\begin{equation}\label{eq:cqr7}
\lambda_i = \lim_{t\to\infty}\frac{1}{t}\psi_t^i, \quad i=1,\ldots,d.
\end{equation}
The difference between the discrete  and the continuous $QR$ method is that
for the first one, the orthonormalization is performed numerically at every discrete time
step, while the continuous $QR$ method maintains the orthogonality via solving differential
equations that encode the orthogonality continuously.

\subsection{Computational Considerations}\label{sec:numi}
In this section, we discuss additional aspects of the computational implementation of discrete and
continuous $QR$ methods to calculate LEs. The application of the discrete $QR$ technique mainly
requires the numerical integration of the SDEs (\ref{eq:SDE3}) and (\ref{eq:D_QR_VE1}). This task
is performed by using standard weak Euler-Maruyama and Milstein schemes, which preserve the
ergodicity property (see \cite{CarBJ10,GroT90,Tal90}).

On the other hand, the numerical integration of the SDEs (\ref{eq:SDE3}), (\ref{eq:cqr3}) and
(\ref{eq:cqr6}) in the computational implementation of the continuous $QR$ technique, must be
performed in such a way that it preserves the orthogonality of the factor $\bfg Q$ in each
integration step. This can be achieved via projected orthogonal schemes which consist of a two-step
process in which first an approximation is computed via any standard scheme, and then the result is
projected into the set of orthogonal matrices \cite{DieRV94}. Again we use the Euler-Maruyama and
Milstein method as in the discrete case.

We have  implemented the two $QR$-methods in \MATLAB. However, to obtain a unique $QR$
factorization in each step, we have modified
%all diagonal elements of $\bfg R$ are
%chosen positive (or alternatively, all diagonal elements are chosen negative), we had to modify
the $QR$ decomposition  provided by \MATLAB to ensure this uniqueness, by forming a diagonal matrix
$\bfg{\mathcal{I}}$ with $\mathcal{I}_{i,i}=\mbox{\rm sign}(R_{i,i})$, for $i=1,\ldots,d$; and then
setting $\bfg Q:=\bfg Q\bfg{\mathcal{I}}$ and $\bfg R:=\bfg{\mathcal{I}}\bfg R$.

%numerical integration of $d(d+1)$ SDEs
%numerical integration of $d(d+2)$ SDEs
%Accuracy and computational time efficiency of the two $QR$ methods with two different solvers is
%studied in the next example.
\subsection{A Simple Numerical Example}\label{sec:exam}
In this subsection we illustrate the described procedures via a simple  strangeness-free SDAE
system in order to compare the computational efficiency and accuracy of both the discrete and
continuous $QR$ methods using the numerical integration schemes Euler-Maruyama and Milstein. The
four numerical methods will be denoted as \emph{D-EM, D-Milstein, C-EM, C-Milstein}, respectively.
The computations are carried out with \MATLAB Version 9.7.0(R2019b) on a computer with CPU Intel
Core i7 composed by 6 cores of 2.20GHz, and 16 GB of RAM.

%Our first example is an SDAE $n=2$ index-1 SDAE example formulated in It\^o sense
%
As a simple example consider the SDAE equation
\begin{equation}\label{eq:ex1}
    \begin{bmatrix}
        1 & 0\\
        0 & 0
    \end{bmatrix}
    \cdot d
    \begin{bmatrix}
        x_1\\
        x_2
    \end{bmatrix}
    =
    \begin{bmatrix}
        -x_2\\
        -\alpha x_1 + \arctan{(x_1)} + x_2
    \end{bmatrix}dt
    +
    \begin{bmatrix}
        (x_1^2+1)^\frac{1}{2}\\
        0
    \end{bmatrix} dw_t,
\end{equation}
with $\alpha\in\mathbb{R}^+$. The nonlinear functions in both the drift and diffusion  part are
continuous on $\mathbb{R}^+$, with continuous and  bounded derivatives, and $w_t$ is a
one-dimensional Wiener process. The underlying SDE of (\ref{eq:ex1}) is
\begin{equation}\label{eq:ex1a}
    d\hat x_t = [-\alpha \hat x_t+\arctan{(\hat x_t)}]dt + (\hat x_t^2+1)^\frac{1}{2}dw_t,
\end{equation}
whose LE exists and can be explicitly represented as the following integral with respect to the
solution of a stationary Fokker–Planck equation (see further details in \cite{Arn98})
\begin{equation}\label{eq:LE-ex1-FP}
    \lambda = -\alpha + \frac{1}{2}\int_{\mathbb{R}}\frac{(\hat x^2-2)}{(\hat x^2+1)}p(x)dx,
\end{equation}
where $p(x)$ is the stationary density of the unique invariant probability law of $\hat x_t$. By
solving numerically (\ref{eq:LE-ex1-FP}) for $\alpha=2$, we obtain the \emph{exact} value of the LE
associated to (\ref{eq:ex1a} and its original SDAE (\ref{eq:ex1}, which is $\lambda=-1.3385$.
%~\cite{GroT90,CarBJ10}.
The accuracy of the $QR$-based methods will be assessed by comparing with  this value as reference.

\begin{figure}[b!]
    \centering
    \includegraphics[scale=0.75]{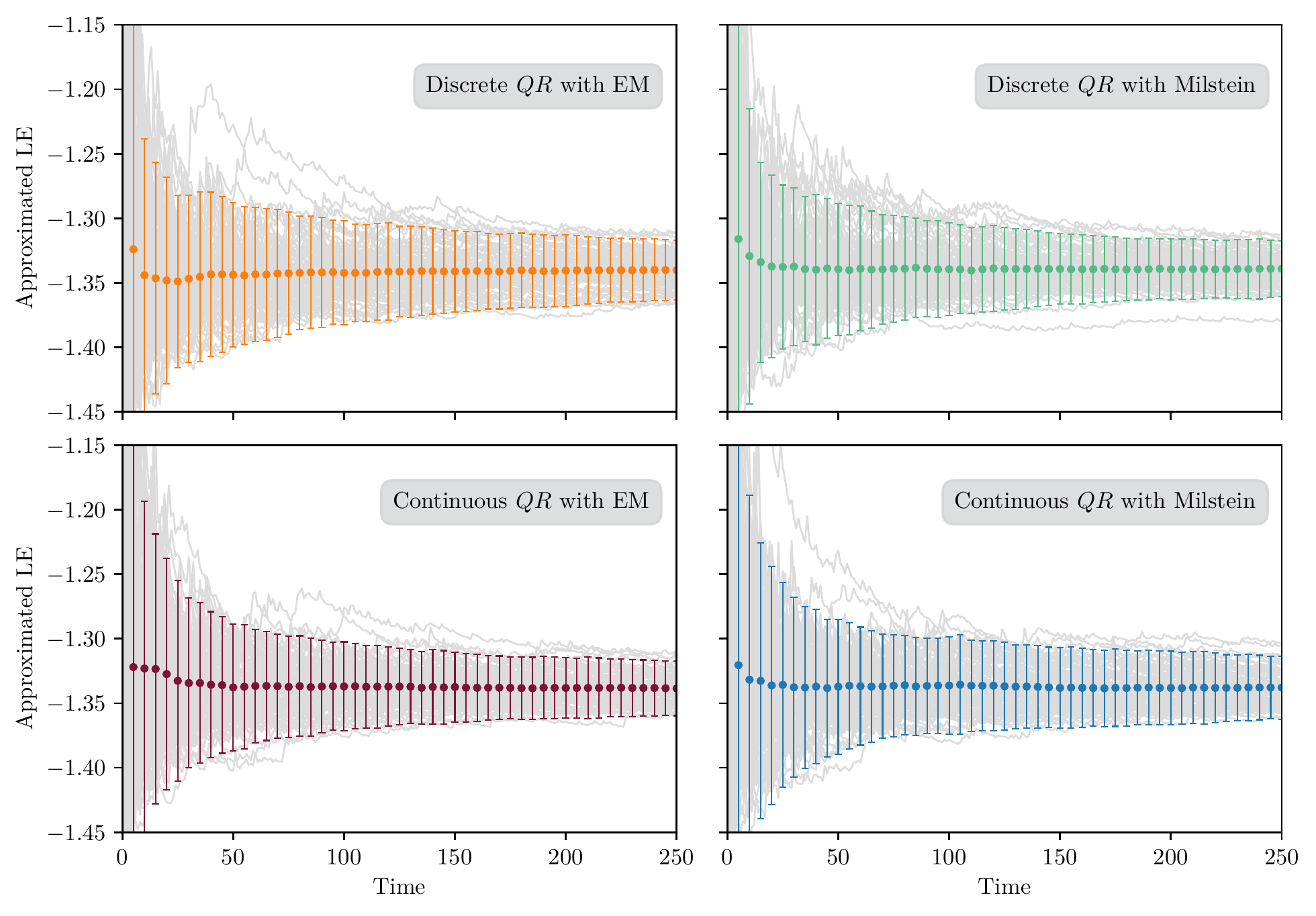}
    \caption{Discrete and continuous $QR$-based approximations of the LE corresponding to
    SDAE (\ref{eq:ex1}) via  Euler-Maruyama and Milstein integrators,
    with a stepsize $h=\exn{1}{-3}$ and $T=250$. The solid circles show mean and the whiskers the
    $95\%$ confidence intervals of the trajectories.}
    \label{fig:LES1}
\end{figure}

\begin{figure}[htb!]
    \centering
    \includegraphics[scale=0.77]{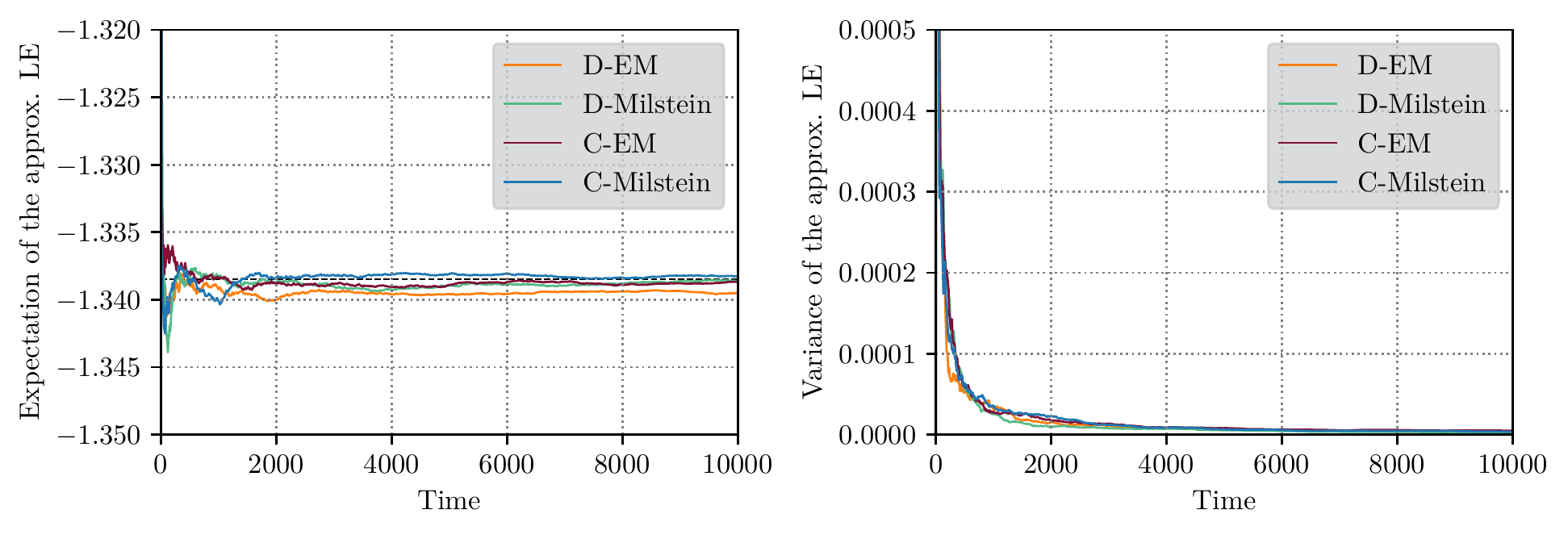}
    \caption{Discrete and continuous $QR$-based approximations of the LE corresponding to  SDAE
    (\ref{eq:ex1}) via  Euler-Maruyama and Milstein integrators, with a stepsize $h=\exn{1}{-3}$
    and $T=10000$. The black dashed line in the left-hand side subplot shows the analytic value of
    $\lambda$.}
    \label{fig:LES2}
\end{figure}

A large number of simulations have been carried out for stepsizes
$h=\exn{1}{-2},\,\exn{9}{-3},\ldots,\,\exn{1}{-3}$ with $T=1000,\,2000,\ldots,\,12000$; to obtain
computed approximations of the LE truncated at the \textit{final time} $t_f:=T$, denoted by
$\lambda_{T}$. To complete our stochastic numerical analysis of the LE, we have calculated the
values of expectation $\mathbb{E}[\lambda_{T}]$, standard deviation $\sigma[\lambda_{T}]$, and
variance $\mathbb{V}[\lambda_{T}]$; estimated from $100$ independent realizations. Some results are
presented in Tables (\ref{tab:e1-d-em}) to (\ref{tab:e1-c-mil}), taking $T=6000,12000,20000$ and
$h=\exn{1}{-1},\,\exn{1}{-2},\,\exn{1}{-3}$.

Observe that the time scale in Figure~\ref{fig:LES1} has been conveniently adjusted to the range
$[0,\,250]$, in order to show the exponential drop of the LE for the different realizations in the
four methods, along the time evolution. While in Figure~\ref{fig:LES2} the time scale has been
adjusted to the range $[0,\,10000]$, to better display the convergence of the mean and  variance
of the LE.

Based on the analytic expression of the LE, given
%in this example by means of the Fokker–Planck
by equation (\ref{eq:LE-ex1-FP}), the LE $\lambda$ can be considered as a deterministic quantity.
According to the numerical results obtained from the four $QR$-based methods, the sequences of
random variables $\lambda_{t_\ell}$ reveal a trend towards null
variance and convergence to the mean as $\ell$ tends to infinity. Such evolution can be seen in
Figure~\ref{fig:LES1}, and more obviously in Figure~\ref{fig:LES2}. For all the methods, an
exponential decay is illustrated in $\mathbb{E}[\lambda_{t_\ell}]$ and
$\mathbb{V}[\lambda_{t_\ell}]$ as $\ell$ tends to infinite. This behavior indicates a mean square
(m.s.) convergence of those sequences to a degenerate random variable, based on the implication
that if $\lambda_{t_\ell} $ is such that $\mathbb{E}[\lambda_{t_\ell}]=\mu_\lambda$,  for all
$\ell$, and $\mathbb{V}[\lambda_{t_\ell}]\xrightarrow[\ell\rightarrow\infty]{}0$, then
$\lambda_{t_\ell} \xrightarrow[\ell\rightarrow\infty]{\textrm{m.s.}} \mu_\lambda$.
This means that the limit of $\lambda_{t_\ell}$ can be interpreted as a deterministic value with
probability $1$. This enables us to state that the stochastic approximations $\lambda_{t_{\ell}}$
converge in m.s. sense to a number (a degenerate random variable), which is expected to represent
the LE $\lambda$.
% very accurate and convergent approximation of the single realizations of $\lambda_{T}$ towards
% the true $\lambda$ can be observed in Figures
% mean of the computed LEs as
% $\overline\lambda_{t_{\ell}}=\frac{1}{100}\sum_{\kappa=1}^{100}\tilde\lambda_{t_{\ell}}^\kappa$
% variance $\lambda_{t_{\ell}}=\frac{1}{99}\sum_{\kappa=1}^{100}[\tilde\lambda_{t_{\ell}}
% ^\kappa-\overline\lambda_{t_{\ell}}]^2$

\begin{figure}[htb!]
    \centering
    \includegraphics[scale=0.7]{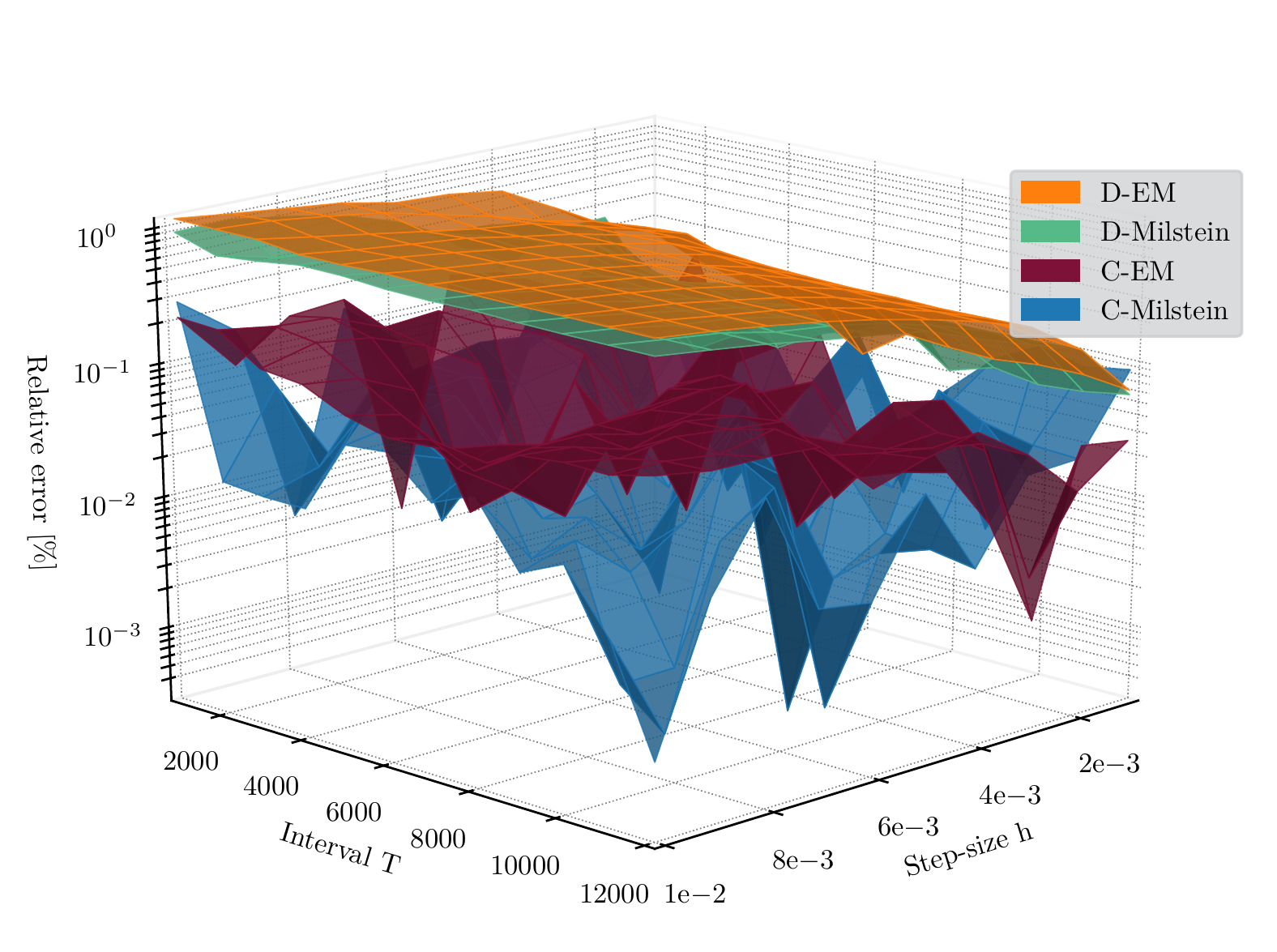}
    \caption{Comparison of relative errors for discrete and continuous $QR$-based approximations of
    the LE corresponding to SDAE (\ref{eq:ex1}) via
    Euler-Maruyama and Milstein integrators, with a range of stepsizes between
    $h=\exn{1}{-2},\ldots,\exn{1}{-3}$; and with $T=1000,\ldots,12000$.}
    \label{fig:Rel_error}
\end{figure}

\begin{figure}[htb!]
    \centering
    \includegraphics[scale=0.7]{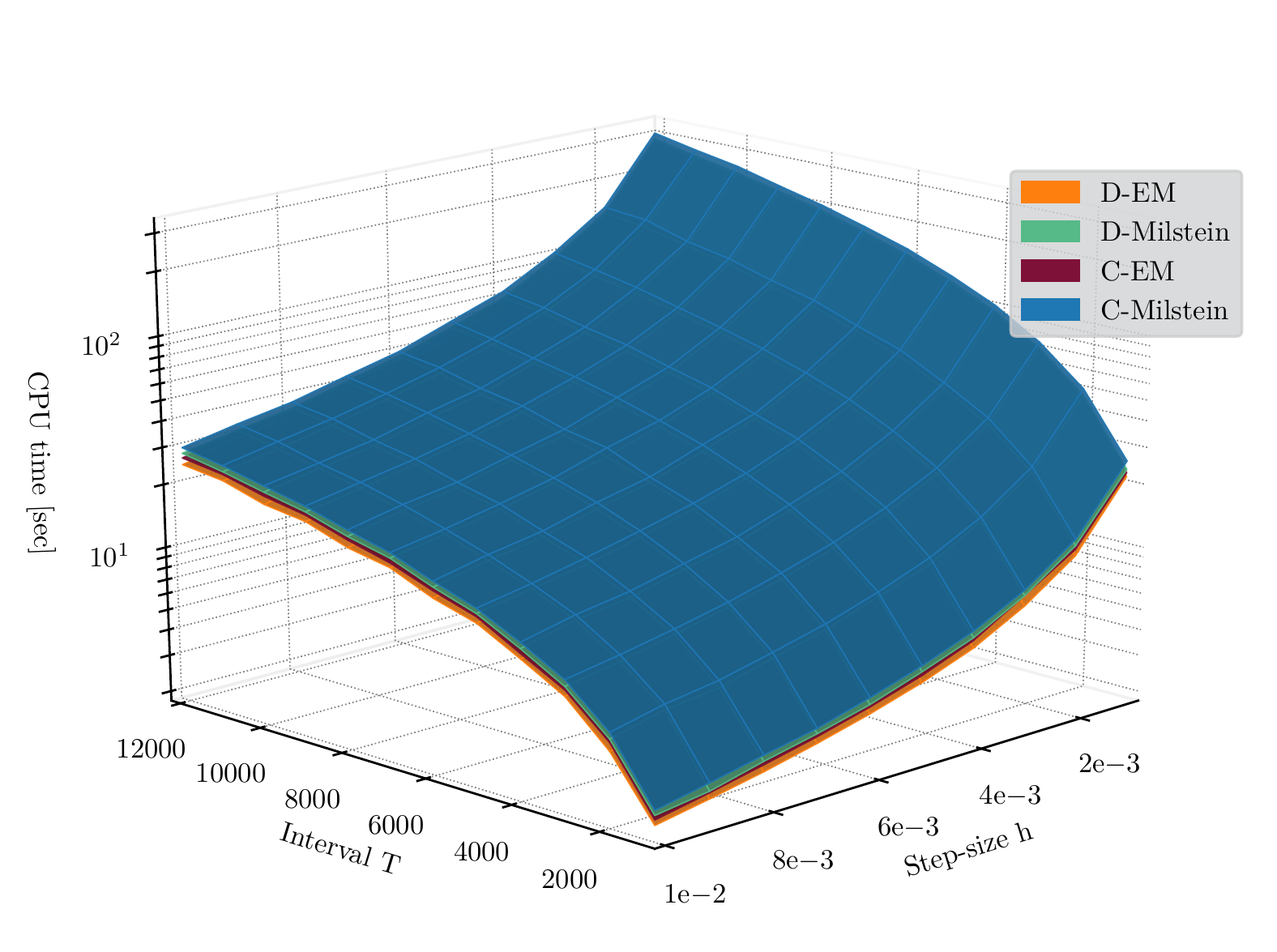}
    \caption{Comparison of the computing-time for discrete and continuous $QR$-based approximations
    of the LE corresponding to SDAE (\ref{eq:ex1}) via
    Euler-Maruyama and Milstein integrators, with a range of stepsizes between
    $h=\exn{1}{-2},\ldots,\exn{1}{-3}$; and with $T=1000,\ldots,12000$.}
    \label{fig:CPU_time}
\end{figure}

In Figure~\ref{fig:Rel_error} we compare the relative error of the accuracy of the four numerical
methods for different stepsize $h$ and time interval $[0,T]$. From this graphical representation,
we observe that continuous methods obtain better results than discrete ones, as expected. We also
observe that the Milstein method has, in general, better accuracy than the Euler-Maruyma scheme,
since its convergence order is higher, but requies more computational time. This latter fact is
evidenced in Figure~\ref{fig:CPU_time}, where a comparison of CPU time (in seconds) is shown for
different values of $h$ and $T$. Here, we observe that all the methods are affected to the same
extent by incrementing the simulation interval $T$, via a logarithmic increment, and by narrowing
the stepsizes $h$, via an exponential increment. A more pronounced difference between the methods
should be observed in higher dimensional systems.

%%%%%%%%%%%%%%%%%%%%%%%%%%%%%%%%%%%%%%%%%%%%%%%%%%%%%%%%%%%
%%%%%%% SECTION 4: STUDY-CASES %%%%%%%%%%%%%%%%%%%%%%%%%%%%
%%%%%%%%%%%%%%%%%%%%%%%%%%%%%%%%%%%%%%%%%%%%%%%%%%%%%%%%%%%
\section{Application of LEs to Power Systems Stability Analysis}\label{sec:apps}
\begin{figure}[htb!]
    \centering
    \includegraphics[scale=0.75]{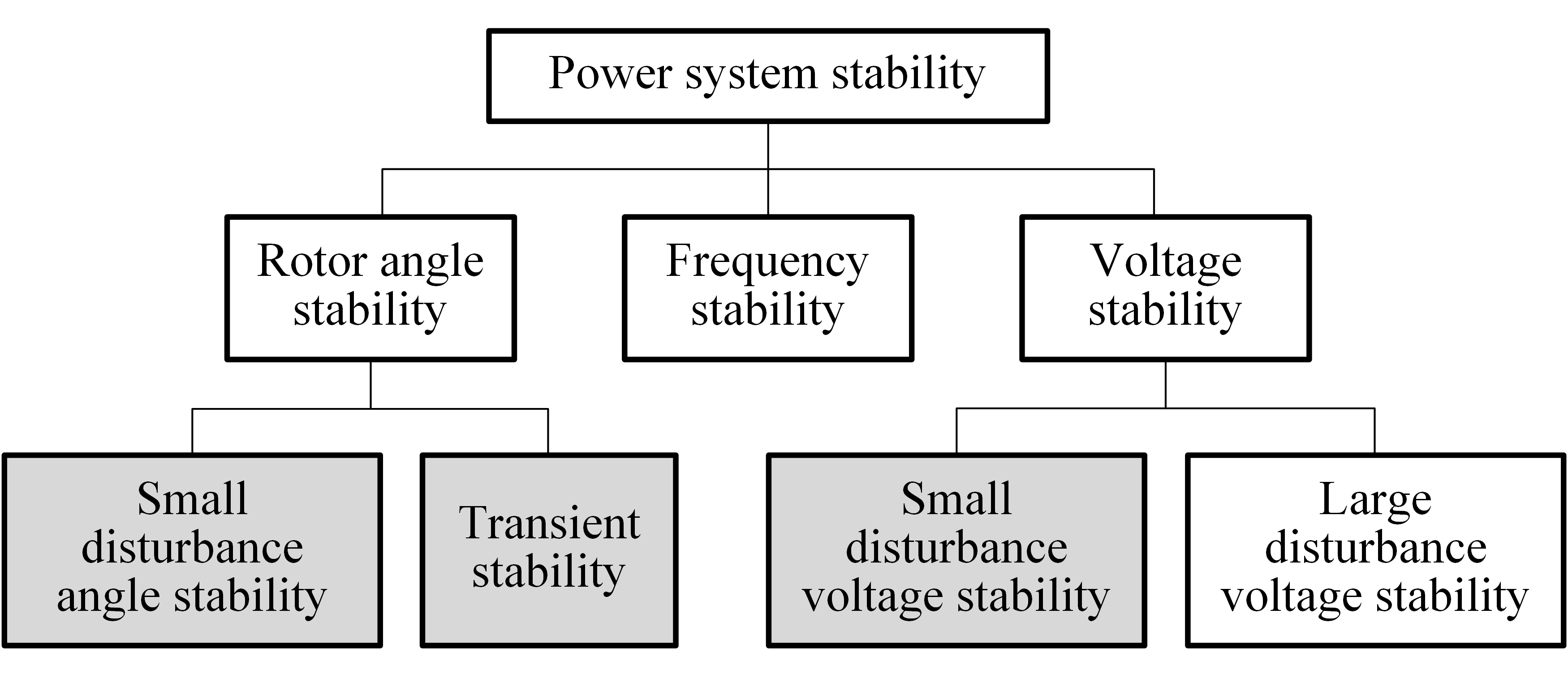}
    \caption{IEEE/CIGRE Power systems stability classification \cite{KunOthers04}.}
    \label{fig:pss_class}
\end{figure}

The concept of stability (based on Lyapunov exponents) in power systems is, in essence, the same as
that for a general dynamical system.  In the literature, power system stability is defined as the
stability to regain an equilibrium state after being subjected to physical disturbances
\cite{Kun94,KunOthers04}. Such equilibrium is characterized through three significant quantities
during the power system operation: angles of nodal voltages, frequency, and nodal voltage
magnitudes. Based on this triplet, there is an entire classification proposed by the
\textit{Institute of Electrical and Electronics Engineers} (IEEE) and the \textit{International
Council on Large Electric Systems} (CIGRE) in \cite{KunOthers04}, which is illustrated in
Figure~\ref{fig:pss_class}.
The test cases presented in this paper are oriented to evaluate the angle and voltage stability of
power systems subjected to small or large disturbances. Studies considering small disturbances are
commonly known as small-signal stability assessment (SSSA). Here, linear stability analysis via
eigenvalues has been one of the traditional analysis tools to predict the degree of stability of
the power system~\cite{Kun94,Sau17}. However, eigenvalue analysis is limited to linear
time-invariant systems or systems close to a stationary solution. When time-varying systems are
tested, as is the case of systems subjected to stochastic disturbances, then eigenvalue analysis is
no longer applicable. On the other hand, the stability analysis of power systems affected by large
disturbances, known as transient stability assessment, is mainly performed with verification
strategies based on time-domain integration \cite{KunOthers04,MacBB11}.

Since the concept of LEs is based on the trajectories of the dynamical systems, the method is an
interesting measure of dynamic stability for power systems under stochastic disturbances in
general. So, testing asymptotic stability of power systems via LEs has become an attractive
approach for the two areas mentioned before, namely, the SSSA of rotor angles and voltages, by
using the linearized set of SDAEs which model the system \cite{VerVK12,VerEKABV16}; and strategies
for the rotor angles via transient analysis using the nonlinear SDAE system and its variational
equation \cite{HayM18}, \cite{WadQA13}. For both cases, asymptotic stability  is checked via
approximations of the largest Lyapunov exponent (LLE) of the system. In particular, a negative LLE
indicates that  the dynamics of the system is asymptotically stable. In the next subsections test
cases are presented that illustrate for strangenss-free SDAE systems the negativity of the LLE.

\subsection{Modeling Power Systems through SDAEs}
Under the assumption of deterministic dynamic behavior, power systems are typically modelled via
a system of quasi-linear DAEs with partitioned variables, see \cite{Kun94},\cite{Sau17}, of the form
\begin{subequations}\label{eq:ps-DAE}
    \begin{align}
        \bfg E_{11}d\bfg x_t^{D_1} &= \bfg f_0^{D_1}(\bfg x_t^{D_1},\bfg x_t^A)dt,
        \label{eq:ps-DAE-d}\\
        0 &= \bfg f_0^A(\bfg x_t^{D_1},\bfg x_t^A),\label{eq:ps-DAE-a}
    \end{align}
\end{subequations}
%
%where $\bfg x_t=[\bfg x_t^{D_1},\bfg x_t^A]^T$, the $\bfg E_1\in\mathbb{R}^{d_1\times n_1}$
%satisfies $\bfg E_1 := [\begin{smallmatrix} \bfg E_{11}& \bfg 0 \end{smallmatrix}]$
where $\bfg E_{11}\in\mathbb{R}^{d_1\times d_1}$ is a diagonal block matrix,
%and $\bfg 0$ is a null matrix of size ${d_1\times a}$,
%inside the leading matrix structure $\bfg E := \bigl[\begin{smallmatrix} \bfg E_{11}&0\\ 0&0
%\end{smallmatrix}\bigr]$,
$\bfg f_0^{D_1}\in\mathcal{C}^1(\mathbb{R}^{d_1+a},\mathbb{R}^{d_1})$,
$\bfg f_0^A\in\mathcal{C}^1(\mathbb{R}^{d_1+a},\mathbb{R}^a)$, $\bfg x_t^{D_1}\in\mathbb{R}^{d_1}$
are the dynamic state variables, and $\bfg x_t^A\in\mathbb{R}^a$ are the algebraic state variables
and we set $n_1=d_1+a$. The DAE system (\ref{eq:ps-DAE}) is  strangeness-free (of index one).

The dynamic behavior of synchronous machines, system controllers, power converters, transmission
lines, or power loads are adequately represented through such a DAE formulation. But in current
real world systems, the dynamic behavior of power systems is affected by disturbances of a
stochastic nature such as renewable stochastic power generation, rotor vibrations in synchronous
machines, stochastic variations of loads, electromagnetic transients, or perturbations originated
by the measurement errors of control devices, see \cite{MilZ13}. Such disturbances can be modeled
through It\^o SDEs of the form
\begin{equation}\label{eq:ps-SDE}
	d\bfg x^{D_2}_t = \bfg f_0^{D_2}(\bfg x^{D_1}_t,\bfg x_t^{D_2},\bfg x_t^A)dt + \bfg
    f_1^{D_2}(\bfg x^{D_1}_t,\bfg x_t^{D_2},\bfg x_t^A)dw_t.
\end{equation}
Here, $\bfg f_0^{D_2}\in\mathcal{C}^1(\mathbb{R}^{d_2+a},\mathbb{R}^{d_2})$ is the drift,
$\bfg f_1^{D_2}\in\mathcal{C}^2(\mathbb{R}^{d_2+a},\mathbb{R}^a)$ is the diffusion, $\bfg
x_t^{D_2}\in\mathbb{R}^{d_2}$ are the stochastic variables, and $w_t$ is the Wiener process. By
combining (\ref{eq:ps-SDE}) and (\ref{eq:ps-DAE}), and assuming that $\bfg x_t^{D_2}$ perturbs
(\ref{eq:ps-DAE-d}) and (\ref{eq:ps-DAE-a}), we obtain a strangeness-free SDAE system of the form
\begin{subequations}\label{eq:ps-SDAE}
    \begin{align}
        \bfg E_{1}d\bfg x_t^{D_1} &= \bfg f_0^{D_1}(\bfg x^{D_1}_t,\bfg x_t^{D_2},\bfg x_t^A)dt,
        \label{eq:ps-SDAE-d}\\
        d\bfg x^{D_2}_t &= \bfg f_0^{D_2}(\bfg x^{D_1}_t,\bfg x_t^{D_2},\bfg x_t^A)dt + \bfg
        f_1^{D_2}(\bfg x^{D_1}_t,\bfg x_t^{D_2},\bfg x_t^A)dw_t,\label{eq:ps-SDAE-s}\\
        0 &= \bfg f_0^A(\bfg x^{D_1}_t,\bfg x_t^{D_2},\bfg x_t^A),\label{eq:ps-SDAE-a}
    \end{align}
\end{subequations}
or in simplified notation  as
\begin{equation}\label{eq:ps-SDAE-comp}
	\bfg E d\bfg x_t = \bfg f_0(\bfg x_t)dt + \bfg f_1(\bfg x_t)dw_t, \quad \bfg
    x_{t_0}=\bfg x_0,
\end{equation}
with
\[
    \bfg E := \left[ \begin{smallmatrix} \bfg E_{11}&0&0\\ 0&\bfg I_{d_2}&0\\0&0&0
    \end{smallmatrix}\right], \ \bfg x_t:=\left [ \begin{array} {c} \bfg x_t^{D_1}\\ \bfg
x^{D_2}_t\\ \bfg x_t^A \end{array} \right ],
\]
drift $\bfg f_0\in\mathcal{C}^1(\mathbb{R}^{n},\mathbb{R}^{n})$,
and  diffusion $\bfg f_1\in\mathcal{C}^2(\mathbb{R}^{n},\mathbb{R}^n)$, where $n=d_1+d_2+a$.

The study-cases presented below are formulated as the form (\ref{eq:ps-SDAE}). An alternative
approach for including the stochastic disturbances is to implement the Wiener process directly in
the underlying ODE of the system, turning them into SDEs (see \cite{CarH17,GeuHL19} for examples).

\subsection{Modeling Stochastic Perturbations}
In this subsection, we discuss  on the modeling  of stochastic variations via SDEs. We employ the
well known mean-reverting process termed \emph{Ornstein–Uhlenbeck} (OU) process
\cite{MilZ13,Gon17}. The SDE which defines the OU process has the form
\begin{equation} \label{eq:OU1}
	d\eta_t =\alpha(\mu-\eta_t)dt + \beta dw_t, \qquad \eta_{t_0}=\eta_0, \qquad t\in\mathbb{I},
\end{equation}
where $\alpha,\mu,\beta\in\mathbb{R}^+$. The OU process is a stationary autocorrelated Gaussian
diffusion process distributed as $\mathcal{N}(\mu,\beta^2/2\alpha)$. Another mean-reverting choice,
similar to the OU process would be the \emph{Cox–Ingersoll–Ross} (CIR) process, whose realizations
are always nonnegative, in fact, it as a sum of squared OU process \cite{All07}.

It is usually recommended to ensure the boundedness of the stochastic variations for the numerical
implementations. In this regard, suitable resources are odd trigonometric functions such as
a $\sin$ or $\arctan$ to guarantee boundedness. For example, if from (\ref{eq:OU1})
we generate a process with a normal distribution $\mathcal{N}(\mu,\sigma^2)$, for $\mu=0$ and
$\sigma^2=0.16$, this value of variance enables us to generate a mean-reverting stochastic
trajectory, whose confidence interval of $95\%$ ($\pm2\sigma$) is inside the threshold of $\pm1$.
Then, through the functions
\begin{equation}\label{eq:OU2}
    \xi(\eta_t)=\sin \eta_t, \qquad \textrm{or} \qquad
    \chi(\eta_t)=\frac{2}{\pi}\arctan\eta_t,
\end{equation}
we obtain a bounded stochastic variation inside the interval $[-1,1]$,
%. According to the system
%(\ref{eq:ps-SDAE}),
and the OU SDEs, that generate the stochastic variations, are represented by (\ref{eq:ps-SDAE-s}).

To couple the parameters of the system in (\ref{eq:ps-SDAE-d}) and (\ref{eq:ps-SDAE-a}) with a
bounded  stochastic disturbance, we use
\[
p(\eta_t) = p_0 + \rho\xi(\eta_t),
\]
where $p_0$ is a constant parameter, $\eta_t$ is the stochastic process that describes the
variations of the parameter, and $\rho\in\mathbb{R}^+$ is a factor that controls the magnitude of
the perturbation.
%We refer to \cite{MilZ13} for further information in this regard.

\subsection{Test Cases}\label{sec:cases}
\begin{figure}[htb]
    \centering
    \includegraphics[width=0.65\textwidth]{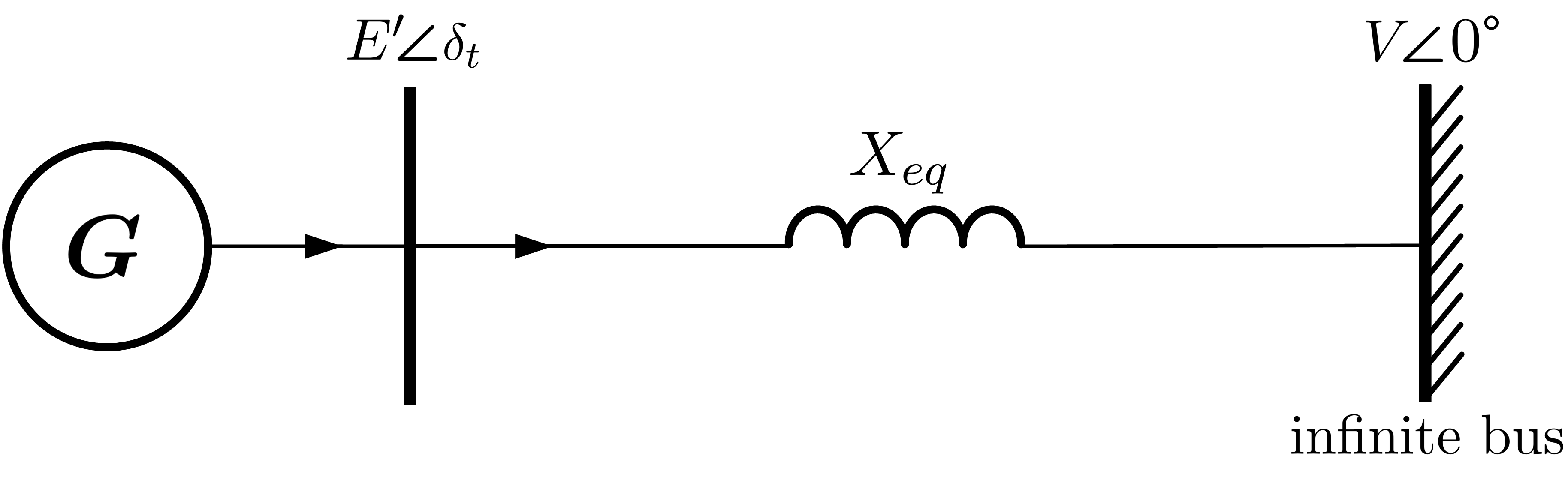}
    \caption{Single-machine infinite-bus (SMIB) scheme.}
    \label{fig:smib}

\end{figure}

In this subsection we present results of our implementation of the $QR$-based methods for the
calculation of LEs at the hand of several test cases of power systems represented by
strangeness-free SDAEs models of so-called single-machine-infinite-bus (SMIB) systems.
% Then, it is possible to estimate the asymptotic
%stability of such power system employing the LLE verification. Specifically, our study-cases make
%use of a representative one-line power system model, called single-machine infinite-bus (SMIB).
This simplified model is frequently used in the area of power systems in order to understand the
local dynamic behavior of a specific machine connected to a complex power network. The SMIB consists
of a synchronous generator connected through a transmission line to a bus with a fixed bus voltage
magnitude and angle, called infinite bus, which represents the grid. A diagram of the system is shown
in Figure~\ref{fig:smib}.

In each test case, we consider a different type of disturbance. For Case 1, the disturbance
is a stochastic load connected to the system. In Case 2, the disturbance is due
to noise caused by a measurement error in a transducer of the machine control system. In both cases,
the maximum disturbance that the system can admit without loosing stability is analyzed, as well
as the effect (positive or negative) of the disturbance for the system in
the stable region. The whole SMIB system, i.e., the synchronous machine, system constraints, and
stochastic disturbances; are modeled by a strangeness-free SDAE system. The dimension of this system
is mainly defined by the type of model used in the synchronous machine; we use a \emph{classical
model} and a \textit{flux-decay model}, see \cite{Kun94,MacBB11,Pot18,Sau17} for detailed
descriptions.

\subsubsection{Case 1: SMIB with stochastic load}
\begin{figure}[htb]
    \centering
    \includegraphics[scale=0.65]{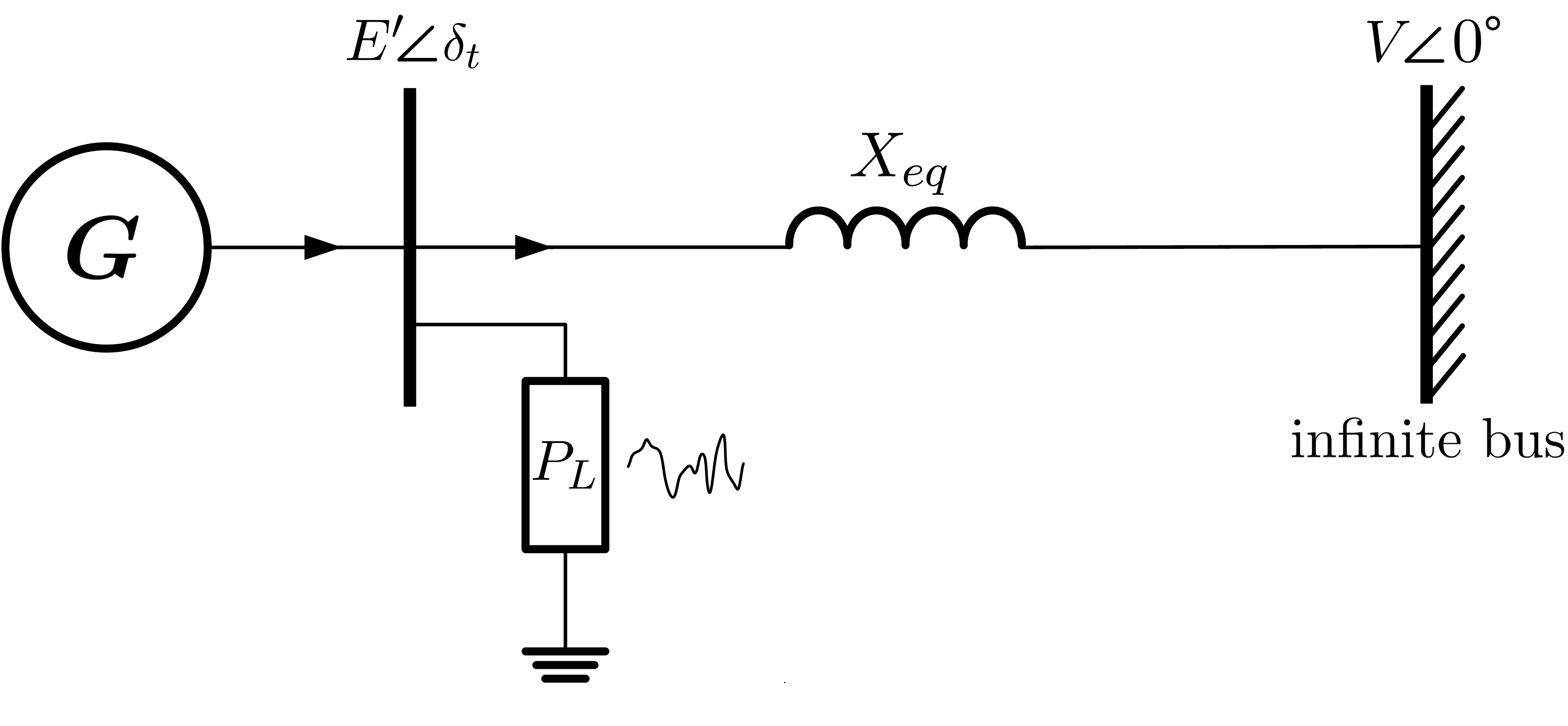}
    \caption{Scheme of a SMIB system with a stochastic load used in Test Case 1.}
    \label{fig:smib-c1}
\end{figure}

In this test case we make use of the LEs to assess the impact of stochastic
disturbances associated with an active power load, over the rotor angle stability of a synchronous
generator. Both the machine and load are connected to the same bus, and this bus in turn, is
linked to the grid through a transmission line. This kind of  SMIB model is typically used to
analyze the effects of renewable energy sources, or aggregated random power consumption, see
Figure~\ref{fig:smib-c1}.
For this version of SMIB system called \emph{classical model}, the dynamic behavior of the
synchronous machine is represented by the swing equations where the rotor angle $\delta_t$ and the
rotor speed $\omega_t$ are the state variables, see \cite{MacBB11,Sau17}. The algebraic constraint
in the system is given by the active power balance, expressed in terms of $P_m$ the mechanical
power, $P_e$ the electrical power, and $P_L$ the constant power consumed by the load. A stochastic
process $\eta_t$ is modeled by an OU SDE. We consider that $\rho\eta_t$ is the stochastic component
of power consumption that perturbs additively the active power balance of the system, where $\rho$
is the size of the disturbance. This leads to the system
%In summary, the system is described by following index-1 SDAE system
\begin{subequations}\label{eq:smib-c1}
    \begin{align}
        d\delta_t &= [\omega_t - \omega_s]dt, \\
        2Hd\omega_t &= [P_m - P_e - K_D(\omega_t-\omega_s)]dt, \\
        \eta_t &= -\alpha \eta_tdt + \beta dw_t, \\
        0 & = \frac{E'V}{X_{eq}}\cos\delta_t + (P_L+\rho\eta_t) - P_e.
    \end{align}
\end{subequations}
By computing the LEs of this SDAE system and checking the LLE, we can determine the maximal
perturbation size $\rho$ (via successive increments of $\rho$) admitted by the SMIB system before
loosing rotor angle stability.  The numerical tests are performed for the values
$P_m=0.8$;
$P_L=0.3$;
$X_{eq}=0.8$;
$H=3.5$;
$K_D=0.4$;
$\omega_s=2\pi50$;
$V=1.0$;
$E'=1.05$;
$\alpha=1.0$;
$\beta=0.4$.
%For the initial values in the state variables: $\delta_t=0.6$; $\omega_t=1.0$.
Most of the values are expressed in the per-unit system (pu) \cite{Pot18}. The
$QR$ methods are executed with  step size $h=\exn{1}{-3}$ and a simulation time $T=20000$.
\begin{figure}[htb!]
    \centering
    \includegraphics[scale=0.65]{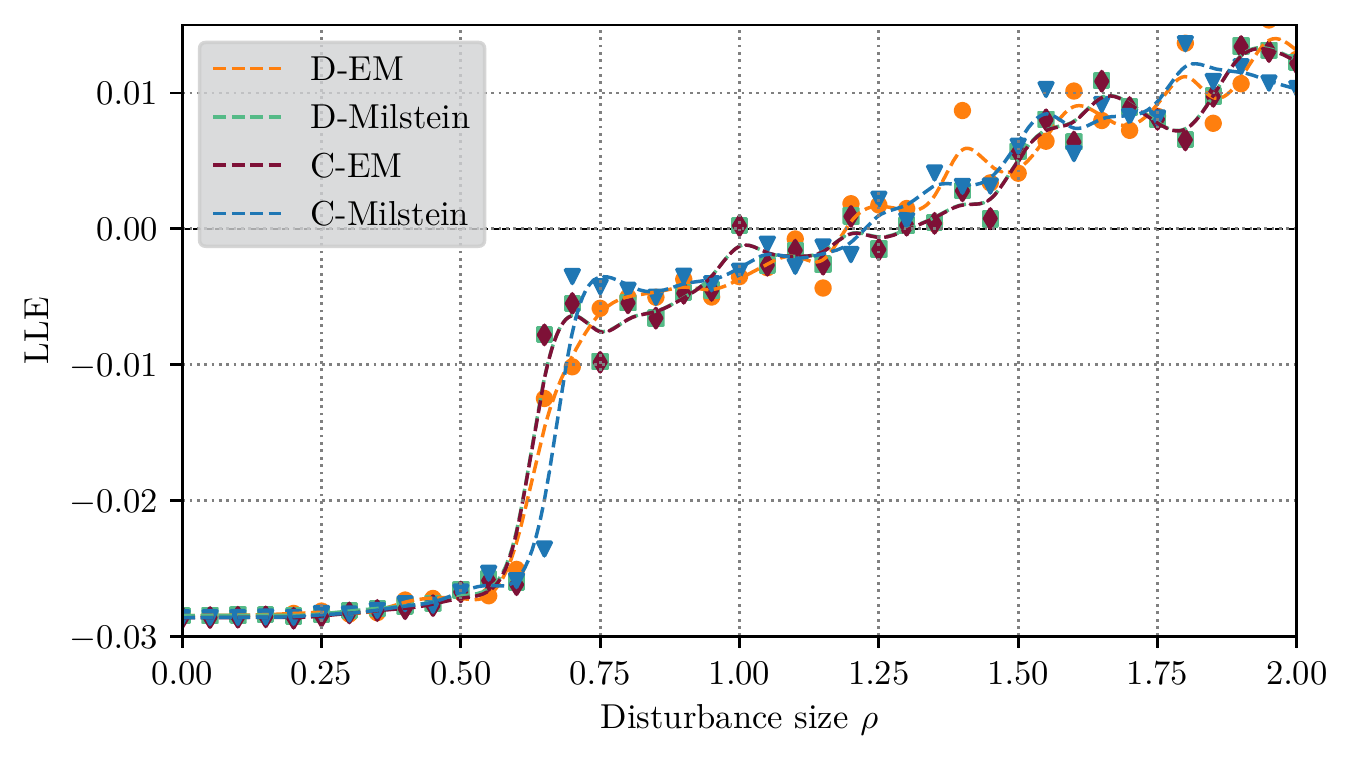}
    \caption{LLE considering different disturbance sizes $\rho$ for the SMIB Test Case 1, tested
    with four $QR$-based methods.}
    \label{fig:LEcase1}
\end{figure}
Figure~\ref{fig:LEcase1} displays the computed LLE utilizing the four $QR$-based methods for
incremental disturbance sizes $\rho=0.00,0.05,\ldots,2.00$. As expected, at $\rho=0.0$ when the
system is not affected by a stochastic disturbance, i.e., the system is deterministic, the computed
LEs match closely with the real parts of the eigenvalues obtained from the Jacobian matrix of the
linearization of (\ref{eq:smib-c1}). When increasing $\rho$, all methods reveal the same
monotonically increasing behavior of the calculated value of the LLE towards the unstable region.
First, there is  a slow increase for $0.00<\rho<0.60$, and then an abrupt increase of the LLE in
the interval $0.60<\rho<0.75$. %, reaching the instability at $\rho\approx0.75$.
% and higher values.
In the interval $0.75<\rho<1.20$, even though the LLE has not yet reached the instability region,
for this particular case the characteristics such as a low damping coefficient and the presence of
the stochastic disturbance, provokes a behavior in the system called \emph{pole slipping}. This is,
in a certain sense, a different kind of instability because the system looses synchronism as it
reaches another equilibrium point near another attractor, see \cite[sec.~5.8]{Sau17} for further
details. The numerical results for this case are presented in Table~\ref{tab:smib_c1}.

This test case shows the large potential of using  LEs as an indicator of instability for nonlinear
power systems. These could be used also in  multi-machine study cases where, however, the
computational complexity has to be reduced, e.g. by model reduction.
% aspect that should be taken
%into account in TSA. %Notwithstanding, its implementation implies the thorough study of further
%concepts, considerations, limitations, etc. that characterize the TSA, and therefore, can be the
%subject of subsequent works in the near future.

\subsubsection{Case 2: SMIB with regulator perturbed by noise}
\begin{figure}[htb]
    \centering
    \includegraphics[scale=0.6]{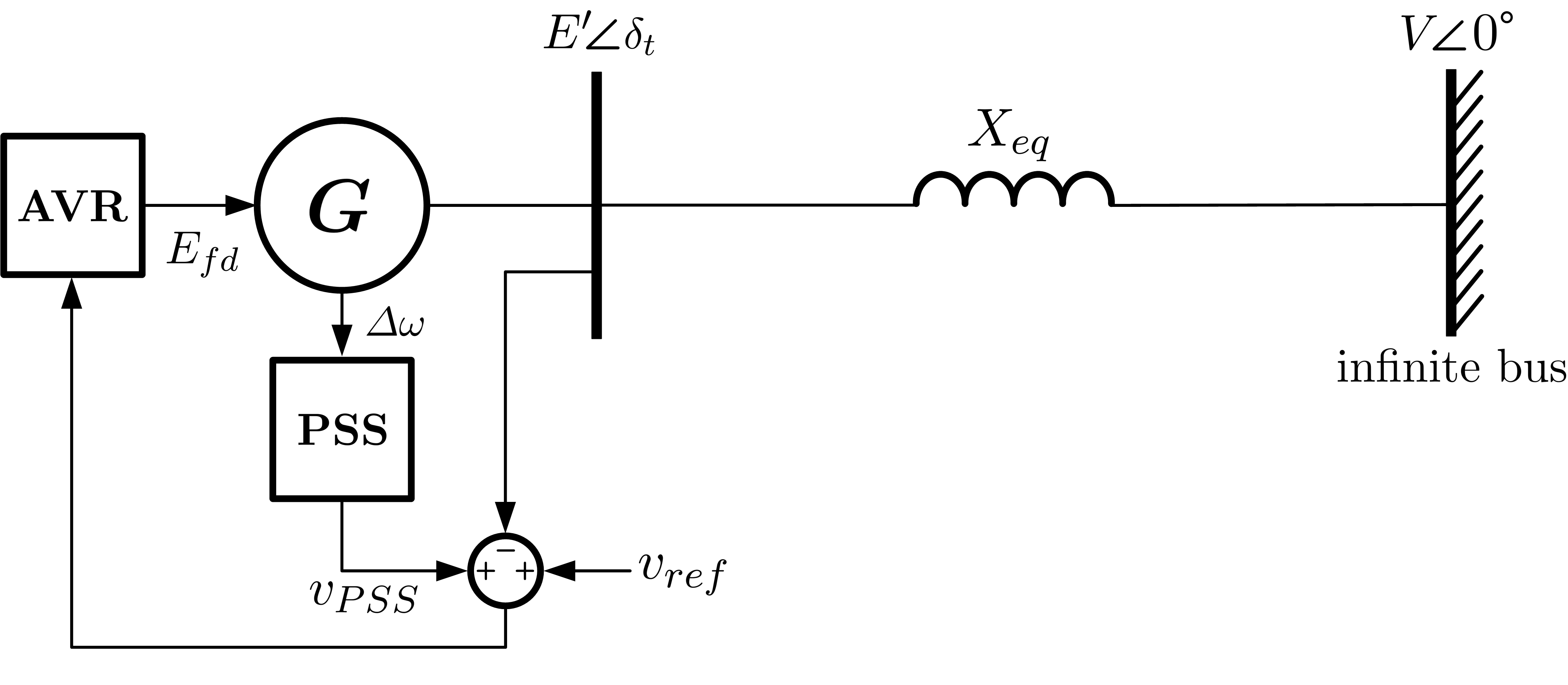}
    \caption{SMIB system scheme equipped with AVR and PSS, corresponding to Test Case 2.}
    \label{fig:smib-c2}
\end{figure}

In this subsection we consider  an SMIB system with a synchronous machine described by a third-order
\emph{flux-decay model}. Here, in addition to the rotor angle $\delta_t$ and the rotor speed
$\omega_t$ associated to the swing equations, the system includes the effect of the field flux
$\psi_{fd}$ described by the field circuit dynamic equations and constraints. In this model the
machine is equipped with an automatic voltage regulator (AVR) to keep the generator output voltage
magnitude in a desirable range, and a power system stabilizer (PSS) to damp out low-frequency
oscillations, see Figure~\ref{fig:smib-c2}.
The AVR and PSS add to the system three more state variables $v_1$, $v_2$, and $v_s$; together
with their corresponding DAEs, which describe the dynamic behavior and constraints of the
controllers into the SMIB system. The resulting model is a nonlinear system of strangeness-free
DAEs. We use the LEs  to analyze the system stability at a specific operation point in the
state-space when it is subjected to small-disturbances. Using the small-signal stability assessment
(SSSA), the set of DAEs that describes the dynamics of the power system is linearized around the
desired operating point. The final result is a linear DAE system. A comprehensive explanation of
this model, its linearization, and reduction to an underlying ODE system can be found in
\cite[ch.~12]{Kun94}. We consider a disturbance of stochastic nature entering in the exciter block
of the AVR as an error of the reference signal \cite{Kun94,VerVK12}, by adding the stochastic
variable $\eta$ to $v_1$ in equation (\ref{eq:noise}). Resolving the algebraic constraints leads to
the linearized system of SDEs
\begin{subequations}\label{eq:smib-c2}
    \begin{align}
        d\Delta\delta &= \omega_s\Delta\omega dt, \\
        2Hd\Delta\omega &= \left[-K_1\Delta\delta-K_D\Delta\omega -K_2\Delta\psi_{fd} + \Delta
        T_m\right]dt, \\ \nonumber
        T_3\Delta\psi_{fd} &= \left[-K_3K_4\Delta\delta-(1+K_3K_6K_A)\Delta\psi_{fd}
        -K_3K_A(1+\rho\eta)\Delta v_1 \right.\\
        & \quad \left. + K_3K_A\Delta v_s\right]dt,\label{eq:noise} \\%+K_3\DeltaE_{fd}
        T_Rd\Delta v_1 &= \left[-K_5\Delta\delta +K_6\Delta\psi_{fd} -\Delta v_1\right]dt, \\
        d\Delta v_2 &= \left[-K_1K_{ST}\Delta\delta -K_DK_{ST}\Delta\omega -K_2K_{ST}\Delta\psi_{fd}
        -\frac{1}{T_W}\Delta v_2 + \frac{K_{ST}}{2H}\Delta T_m\right]dt, \\ \nonumber
        T_2 d\Delta v_s &= \left[-K_1K_{ST}T_1\Delta\delta - K_DK_{ST}T_1\Delta\omega
        -K_2K_{ST}T_1\Delta\psi_{fd}+ \left(\frac{T_1}{T_W}+1 \right)\Delta v_2 \right.\\
        & \quad \left. -\frac{1}{T_2}\Delta v_s+\frac{K_{ST}T_1}{2H}\Delta T_m\right]dt, \\
        d\eta &= -\alpha\eta dt + \beta dw,
    \end{align}
\end{subequations}
where $\Delta\delta, \ \Delta\omega, \ \Delta\psi_{fd}, \ \Delta v_1, \ \Delta v_2, \ \Delta v_s,
\ \eta$ are the state variables of the linear underlying SDE system (for simplicity, the
subscript $t$ has been omitted in this formulation). Once again, the stochastic perturbation
is generated via an OU SDE, and the size of the perturbation is controlled by the parameter $\rho$.
The numerical analysis is done for the values
$\omega_s=2\pi60$;
%$V=1.0$
$H=3.0$;
%$X_{eq}=0.8$;
$K_1=1.591$;
$K_2=1.50$;
$K_D=0.0$;
$K_3=0.333$;
$K_4=1.8$;
$K_5=-0.12$;
$K_6=0.3$;
$K_A=200.0$;
$T_R=0.02$;
$K_{ST}=9.5$;
$T_1=0.154$;
$T_2=0.033$;
$T_3=1.91$;
$T_W=1.4$;
$\alpha=1.0$;
$\beta=0.4$;
$\Delta T_m=0.0$.

\begin{figure}[htb!]
    \centering
    \includegraphics[scale=0.75]{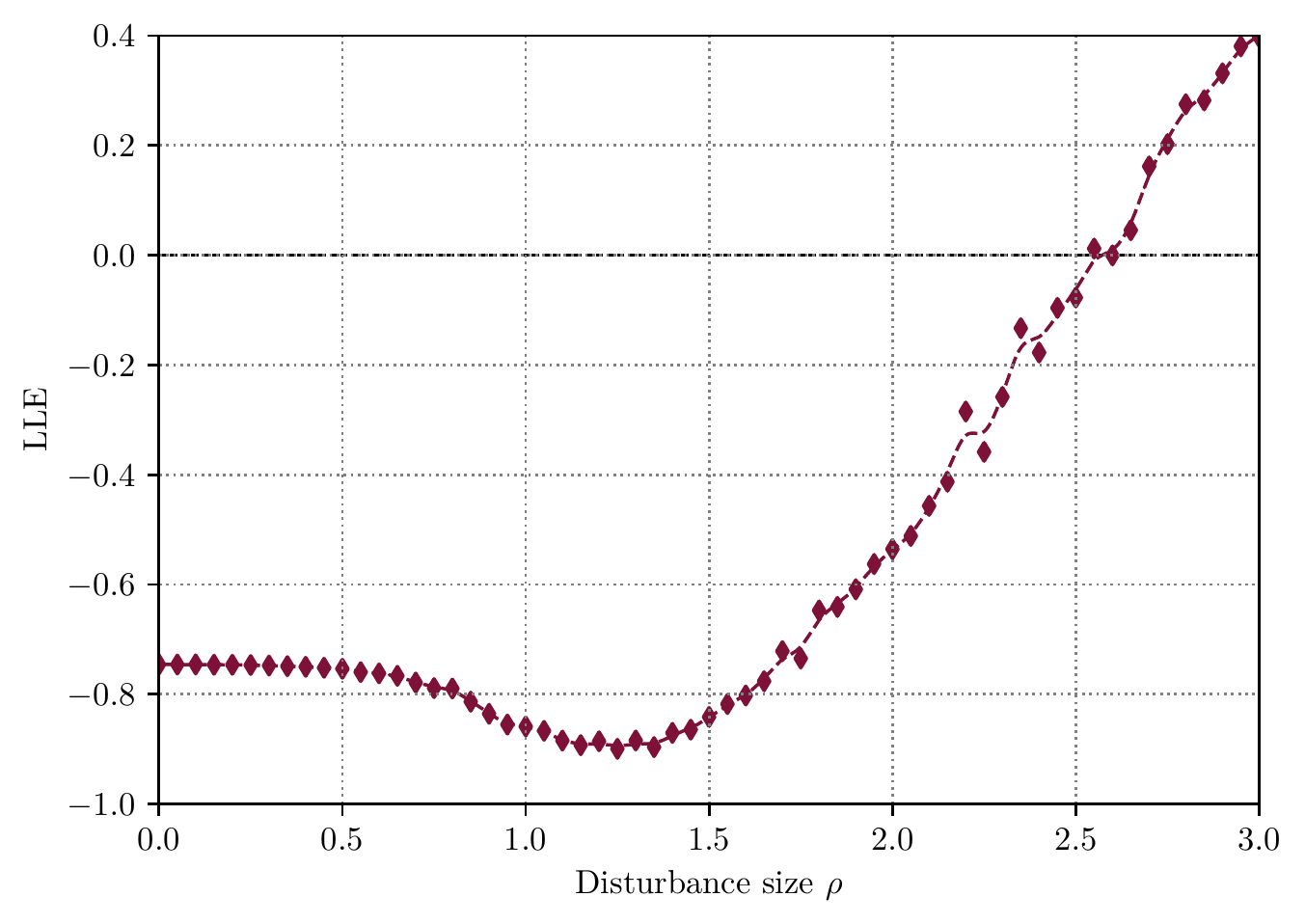}
    \caption{Computed LLE for the dimension 7 SMIB system of Test Case 2, considering different
    disturbance sizes and using the continuous Euler-Maruyama $QR$ method.}
    \label{fig:LEcase2}
\end{figure}

Based on the analysis of Section (\ref{sec:exam}), we only consider the continuous  Euler-Maruyama
QR method. The results of computing the LLE of the SMIB system for incremental values of the
perturbation size $\rho$, are presented graphically in Figure~\ref{fig:LEcase2}.
The values of the LLE when increasing perturbation size $\rho$ clearly mark
four defined intervals. In the leftmost interval with $0.00<\rho<0.40$, the calculated LLE is
practically constant and equal to the real part of the right-most eigenvalue from the deterministic
system. In this region, there is no impact of the disturbance on the system stability. In the
interval $0.40<\rho<1.30$ a curious situation occurs, as the size of the disturbance
increases, the distance from the LLE to the positive region increases, in other words, the noise
improves the stability of the system. In the interval $1.30<\rho<2.60$, the situation changes
completely, and the LLE converges to zero.
Finally, from $\rho\approx2.60$ onwards, the system is unstable. Table \ref{tab:smib_c2} shows
the numerical values of this test case.

\begin{figure}[htb!]
    \centering
    \includegraphics[scale=0.75]{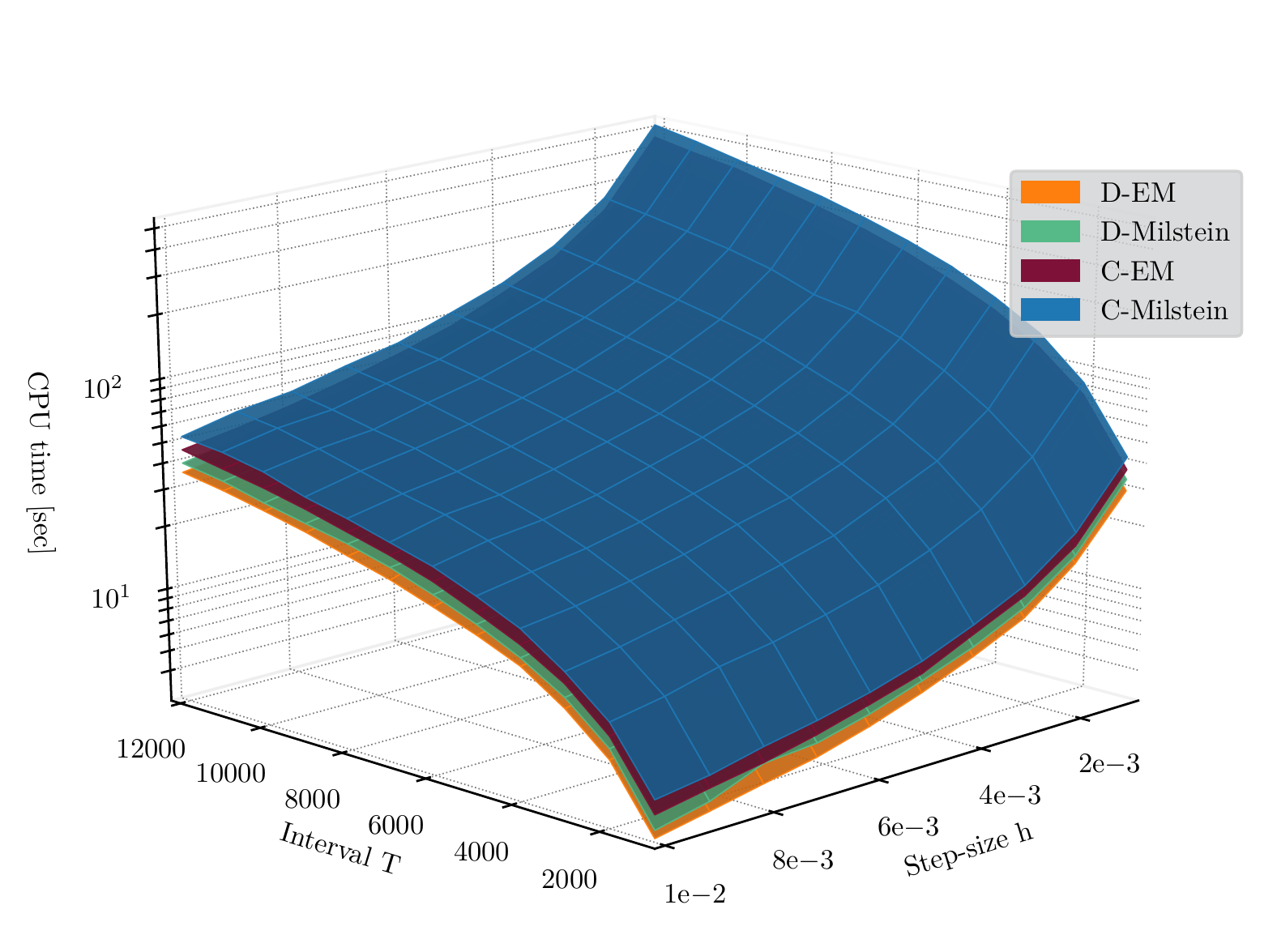}
    \caption{Computing-time comparison of LE calculation for the dimension 7 SMIB Test Case 2.
    Comparison performed for the four $QR$ methods in a range of step sizes between
    $h=[\exn{1}{-2},\exn{1}{-3}]$ and with $T=[1000,12000]$.}
    \label{fig:LE_time_c2}
\end{figure}

Finally, we have evaluated the computing-times for this 7-dimensional test case. The results are
shown in Figure~\ref{fig:LE_time_c2}. Although the computational cost for all method is similar for
the different methods as a factor of the step sizes $h$ and time interval $[0,T]$, the
computational costs strongly increase.
%In fact, while the discrete methods solve $d(d+1)$ SDEs per integration-step, in
%continuous methods $d(d+2)$ SDEs are solved.
%

%%%%%%%%%%%%%%%%%%%%%%%%%%%%%%%%%%%%%%%%%%%%%%%%%%%%%%%%%%%
%%%%%%% SECTION 5: CONCLUSIONS %%%%%%%%%%%%%%%%%%%%%%%%%%%%
%%%%%%%%%%%%%%%%%%%%%%%%%%%%%%%%%%%%%%%%%%%%%%%%%%%%%%%%%%%
\section{Conclusions}\label{sec:conclu}
We have revisited the theory of strangeness-free SDAE systems, as well as the concepts of LEs
associated with the RDEs generated via such SDAEs. We have adapted and implemented stochastic
versions of continuous and discrete $QR$-based methods to calculate approximations of the LEs, and
assessed them by using Euler-Maruyama and Milstein schemes over the corresponding underlying SDE.
The results obtained from our numerical experiments illustrate the approximations of the
corresponding LE converge to degenerate random variables, i.e., the LE can be interpreted as a
deterministic value, since in the limit the variance of the approximations tends to zero. Both
$QR$-based method provide reliable results, but in general, continuous methods provide better
accuracy than the discrete counterpart at the expenses of higher computational cost and higher
memory requirement. We have illustrated the $QR$-based methods for SMIB power system problems and
shown the usefulness of the LEs as a stability indicator for the rotor angle and voltage stability
analysis of power systems affected by bounded stochastic disturbances.
%The numerical experiments evidence this
%assertion considering different intensities for the random fluctuations.

As future work, we suggest the use of discretization schemes for SDAEs in order to directly apply
the numerical integration to the SDAE system. Furthermore, methods for computing the LEs based
on \textit{Singular Value Decompositions}, a combination with model reduction, and a careful
comparison with QR-based methods would be of interest. Concerning to the applications to power
systems and dynamical network systems in general, stability assessment of large-scale cases are
remarkable works to be performed.
%try other stability index such as Bohl exponents
%For more details, the readers are suggested to refer [??].

\section*{Acknowledgment}
A. Gonz\'{a}lez-Zumba acknowledges the support of Secretar\'{i}a Nacional de Ciencia y
Tecnolog\'{i}a SENESCYT (Ecuador), through the scholarship \textit{``Becas de Fomento al
Talento Humano''}, and \emph{Deutsche Forschungsgemeinschaft} through Collaborative Research Centre
Transregio. SFB TRR 154.
P. Fern\'andez-de-C\'ordoba was partially supported by grant no. RTI2018-102256-B-I00 (Spain).
J.-C. Cort\'{e}s acknowledges the support by the Spanish Ministerio de Econom\'{i}a, Industria y
Competitividad (MINECO), the Agencia Estatal de Investigaci\'{o}n (AEI), and Fondo Europeo de
Desarrollo Regional (FEDER UE) grant MTM2017--89664--P.
V. Mehrmann was partially supported by \emph{Deutsche Forschungsgemeinschaft} through the
Excellence Cluster {Math$^+$} in Berlin, and Priority Program 1984 ``Hybride und multimodale
Energiesysteme:  Systemtheoretische Methoden für die Transformation und den Betrieb komplexer
Netze''.

%\appendix
\section*{Appendices}
In this Appendix we present the numerical values for several different simulations.
\renewcommand{\thetable}{A.\arabic{table}}%
\renewcommand{\thefigure}{A.\arabic{figure}}%

\begin{table}[htb!]
    \centering
    \footnotesize
    \begin{tabular}{rrrrrrr}
    \toprule[1.0pt]
    \multicolumn{1}{c}{$T$} & \multicolumn{1}{c}{$h$} & \multicolumn{1}{c}{$\mathbb{E}[\lambda_T]$}
    &   \multicolumn{1}{c}{$\sigma[\lambda_T]$} &   \multicolumn{1}{c}{$\mathbb{V}[\lambda_{T}]$} &
    \multicolumn{1}{c}{Rel. error [\%]} &   \multicolumn{1}{c}{CPU-time [sec]}\\
    \toprule[1.0pt]
    $6000$	&	$\exn{1}{-1}$	&	$-1.51906$	&	$0.00476$	&	$\exn{2.266}{-5}$	&	$13.48994$	&	$0.9827$	\\
$6000$	&	$\exn{1}{-2}$	&	$-1.35231$	&	$0.00296$	&	$\exn{8.734}{-6}$	&	$1.03179$	&	$12.6422$	\\
$6000$	&	$\exn{1}{-3}$	&	$-1.33874$	&	$0.00234$	&	$\exn{5.483}{-6}$	&	$0.01807$	&	$125.6875$	\\
\midrule[0.3pt]
$12000$	&	$\exn{1}{-1}$	&	$-1.51870$	&	$0.00334$	&	$\exn{1.116}{-5}$	&	$13.46268$	&	$2.0336$	\\
$12000$	&	$\exn{1}{-2}$	&	$-1.35217$	&	$0.00184$	&	$\exn{3.392}{-6}$	&	$1.02160$	&	$24.9517$	\\
$12000$	&	$\exn{1}{-3}$	&	$-1.33920$	&	$0.00161$	&	$\exn{2.579}{-6}$	&	$0.05206$	&	$252.4366$	\\
\midrule[0.3pt]
$20000$	&	$\exn{1}{-1}$	&	$-1.51780$	&	$0.00262$	&	$\exn{6.858}{-6}$	&	$13.39551$	&	$3.3806$	\\
$20000$	&	$\exn{1}{-2}$	&	$-1.35236$	&	$0.00139$	&	$\exn{1.944}{-6}$	&	$1.03533$	&	$41.2975$	\\
$20000$	&	$\exn{1}{-3}$	&	$-1.33936$	&	$0.00133$	&	$\exn{1.781}{-6}$	&	$0.06437$	&	$416.3228$	\\

    \toprule[1.0pt]
    \end{tabular}
    \caption{Numerical results of the calculated LE for SDAE system (\ref{eq:ex1}) computed via
    Discrete $QR$-EM method.}
    \label{tab:e1-d-em}
\end{table}

\begin{table}[htb!]
    \centering
    \footnotesize
    \begin{tabular}{rrrrrrr}
    \toprule[1.0pt]
    \multicolumn{1}{c}{$T$} & \multicolumn{1}{c}{$h$} & \multicolumn{1}{c}{$\mathbb{E}[\lambda_T]$}
    &   \multicolumn{1}{c}{$\sigma[\lambda_T]$} &   \multicolumn{1}{c}{$\mathbb{V}[\lambda_{T}]$} &
    \multicolumn{1}{c}{Rel. error [\%]} &   \multicolumn{1}{c}{CPU-time [sec]}\\
    \toprule[1.0pt]
    $6000$	&	$\exn{1}{-1}$	&	$-1.47000$	&	$0.00356$	&	$\exn{1.267}{-5}$	&	$9.82471$	&	$1.3947$	\\
$6000$	&	$\exn{1}{-2}$	&	$-1.34911$	&	$0.00205$	&	$\exn{4.217}{-6}$	&	$0.79239$	&	$14.2033$	\\
$6000$	&	$\exn{1}{-3}$	&	$-1.33883$	&	$0.00185$	&	$\exn{3.415}{-6}$	&	$0.02480$	&	$139.4657$	\\
\midrule[0.3pt]
$12000$	&	$\exn{1}{-1}$	&	$-1.46914$	&	$0.00249$	&	$\exn{6.202}{-6}$	&	$9.75996$	&	$2.8596$	\\
$12000$	&	$\exn{1}{-2}$	&	$-1.34925$	&	$0.00186$	&	$\exn{3.466}{-6}$	&	$0.80302$	&	$27.8426$	\\
$12000$	&	$\exn{1}{-3}$	&	$-1.33889$	&	$0.00176$	&	$\exn{3.093}{-6}$	&	$0.02924$	&	$280.9451$	\\
\midrule[0.3pt]
$20000$	&	$\exn{1}{-1}$	&	$-1.46973$	&	$0.00159$	&	$\exn{2.517}{-6}$	&	$9.80448$	&	$4.6544$	\\
$20000$	&	$\exn{1}{-2}$	&	$-1.34924$	&	$0.00141$	&	$\exn{1.982}{-6}$	&	$0.80274$	&	$46.9282$	\\
$20000$	&	$\exn{1}{-3}$	&	$-1.33915$	&	$0.00130$	&	$\exn{1.699}{-6}$	&	$0.04854$	&	$465.2272$	\\

    \toprule[1.0pt]
    \end{tabular}
    \caption{Numerical results of the calculated LE for SDAE system (\ref{eq:ex1}) computed via
    Discrete $QR$-Milstein method.}
    \label{tab:e1-d-mil}
\end{table}

\begin{table}[htb!]
    \centering
    \footnotesize
    \begin{tabular}{rrrrrrr}
    \toprule[1.0pt]
    \multicolumn{1}{c}{$T$} & \multicolumn{1}{c}{$h$} & \multicolumn{1}{c}{$\mathbb{E}[\lambda_T]$}
    &   \multicolumn{1}{c}{$\sigma[\lambda_T]$} &   \multicolumn{1}{c}{$\mathbb{V}[\lambda_{T}]$} &
    \multicolumn{1}{c}{Rel. error [\%]} &   \multicolumn{1}{c}{CPU-time [sec]}\\
    \toprule[1.0pt]
    $6000$	&	$\exn{1}{-1}$	&	$-1.35864$	&	$0.00326$	&	$\exn{1.061}{-5}$	&	$1.50459$	&	$1.3841$	\\
$6000$	&	$\exn{1}{-2}$	&	$-1.34005$	&	$0.00278$	&	$\exn{7.743}{-6}$	&	$0.11616$	&	$13.6211$	\\
$6000$	&	$\exn{1}{-3}$	&	$-1.33822$	&	$0.00277$	&	$\exn{7.651}{-6}$	&	$0.02128$	&	$135.0360$	\\
\midrule[0.3pt]
$12000$	&	$\exn{1}{-1}$	&	$-1.35932$	&	$0.00226$	&	$\exn{5.091}{-6}$	&	$1.55512$	&	$2.7334$	\\
$12000$	&	$\exn{1}{-2}$	&	$-1.33999$	&	$0.00186$	&	$\exn{3.459}{-6}$	&	$0.11134$	&	$26.9335$	\\
$12000$	&	$\exn{1}{-3}$	&	$-1.33813$	&	$0.00159$	&	$\exn{2.535}{-6}$	&	$0.02786$	&	$272.8842$	\\
\midrule[0.3pt]
$20000$	&	$\exn{1}{-1}$	&	$-1.35888$	&	$0.00196$	&	$\exn{3.835}{-6}$	&	$1.52252$	&	$4.4753$	\\
$20000$	&	$\exn{1}{-2}$	&	$-1.34010$	&	$0.00096$	&	$\exn{9.306}{-7}$	&	$0.11965$	&	$45.0326$	\\
$20000$	&	$\exn{1}{-3}$	&	$-1.33807$	&	$0.00148$	&	$\exn{2.187}{-6}$	&	$0.03224$	&	$465.5119$	\\

    \toprule[1.0pt]
    \end{tabular}
    \caption{Numerical results of the calculated LE for SDAE system (\ref{eq:ex1}) computed via
    Continuous $QR$-EM method.}
    \label{tab:e1-c-em}
\end{table}

\begin{table}[htb!]
    \centering
    \footnotesize
    \begin{tabular}{rrrrrrr}
    \toprule[1.0pt]
    \multicolumn{1}{c}{$T$} & \multicolumn{1}{c}{$h$} & \multicolumn{1}{c}{$\mathbb{E}[\lambda_T]$}
    &   \multicolumn{1}{c}{$\sigma[\lambda_T]$} &   \multicolumn{1}{c}{$\mathbb{V}[\lambda_{T}]$} &
    \multicolumn{1}{c}{Rel. error [\%]} &   \multicolumn{1}{c}{CPU-time [sec]}\\
    \toprule[1.0pt]
    $6000$	&	$\exn{1}{-1}$	&	$-1.33950$	&	$0.00287$	&	$\exn{8.228}{-6}$	&	$0.07468$	&	$1.5009$	\\
$6000$	&	$\exn{1}{-2}$	&	$-1.33767$	&	$0.00253$	&	$\exn{6.386}{-6}$	&	$0.06235$	&	$14.8586$	\\
$6000$	&	$\exn{1}{-3}$	&	$-1.33810$	&	$0.00232$	&	$\exn{5.402}{-6}$	&	$0.02998$	&	$147.5276$	\\
\midrule[0.3pt]
$12000$	&	$\exn{1}{-1}$	&	$-1.33931$	&	$0.00259$	&	$\exn{6.692}{-6}$	&	$0.06075$	&	$3.0427$	\\
$12000$	&	$\exn{1}{-2}$	&	$-1.33864$	&	$0.00121$	&	$\exn{1.460}{-6}$	&	$0.01053$	&	$29.3666$	\\
$12000$	&	$\exn{1}{-3}$	&	$-1.33769$	&	$0.00153$	&	$\exn{2.329}{-6}$	&	$0.06087$	&	$299.2588$	\\
\midrule[0.3pt]
$20000$	&	$\exn{1}{-1}$	&	$-1.33990$	&	$0.00182$	&	$\exn{3.296}{-6}$	&	$0.10453$	&	$4.9331$	\\
$20000$	&	$\exn{1}{-2}$	&	$-1.33828$	&	$0.00152$	&	$\exn{2.310}{-6}$	&	$0.01654$	&	$49.5662$	\\
$20000$	&	$\exn{1}{-3}$	&	$-1.33853$	&	$0.00140$	&	$\exn{1.960}{-6}$	&	$0.00258$	&	$505.4304$	\\

    \toprule[1.0pt]
    \end{tabular}
    \caption{Numerical results of the calculated LE for SDAE system (\ref{eq:ex1}) computed via
    Continuous $QR$-Milstein method.}
    \label{tab:e1-c-mil}
\end{table}

\begin{table}[htb!]
    \centering
    \footnotesize
    \begin{tabular}{rrrrrrrrrr}
    \toprule[1.0pt]
    \multicolumn{1}{c}{$\rho$} & \multicolumn{1}{c}{D-EM} & \multicolumn{1}{c}{D-Mil} &
    \multicolumn{1}{c}{C-EM} & \multicolumn{1}{c}{C-Mil} &
    \multicolumn{1}{c}{$\rho$} & \multicolumn{1}{c}{D-EM} & \multicolumn{1}{c}{D-Mil} &
    \multicolumn{1}{c}{C-EM} & \multicolumn{1}{c}{C-Mil}\\
    \toprule[1.0pt]
    $0.00$	&	$-0.02849$	&	$-0.02849$	&	$-0.02864$	&	$-0.02864$	&	$1.05$	&	$-0.00287$	&	$-0.00267$	&	$-0.00268$	&	$-0.00114$	\\
$0.05$	&	$-0.02848$	&	$-0.02847$	&	$-0.02863$	&	$-0.02863$	&	$1.10$	&	$-0.00075$	&	$-0.00159$	&	$-0.00158$	&	$-0.00278$	\\
$0.10$	&	$-0.02843$	&	$-0.02845$	&	$-0.02860$	&	$-0.02863$	&	$1.15$	&	$-0.00436$	&	$-0.00259$	&	$-0.00261$	&	$-0.00136$	\\
$0.15$	&	$-0.02845$	&	$-0.02841$	&	$-0.02857$	&	$-0.02859$	&	$1.20$	&	$0.00184$	&	$0.00093$	&	$0.00091$	&	$-0.00191$	\\
$0.20$	&	$-0.02832$	&	$-0.02852$	&	$-0.02867$	&	$-0.02863$	&	$1.25$	&	$0.00175$	&	$-0.00149$	&	$-0.00154$	&	$0.00217$	\\
$0.25$	&	$-0.02815$	&	$-0.02836$	&	$-0.02852$	&	$-0.02833$	&	$1.30$	&	$0.00149$	&	$0.00029$	&	$0.00026$	&	$0.00059$	\\
$0.30$	&	$-0.02837$	&	$-0.02811$	&	$-0.02828$	&	$-0.02837$	&	$1.35$	&	$0.00038$	&	$0.00044$	&	$0.00040$	&	$0.00409$	\\
$0.35$	&	$-0.02827$	&	$-0.02795$	&	$-0.02811$	&	$-0.02811$	&	$1.40$	&	$0.00870$	&	$0.00284$	&	$0.00279$	&	$0.00311$	\\
$0.40$	&	$-0.02734$	&	$-0.02778$	&	$-0.02797$	&	$-0.02762$	&	$1.45$	&	$0.00338$	&	$0.00072$	&	$0.00075$	&	$0.00314$	\\
$0.45$	&	$-0.02722$	&	$-0.02758$	&	$-0.02775$	&	$-0.02773$	&	$1.50$	&	$0.00409$	&	$0.00570$	&	$0.00564$	&	$0.00607$	\\
$0.50$	&	$-0.02676$	&	$-0.02658$	&	$-0.02674$	&	$-0.02675$	&	$1.55$	&	$0.00644$	&	$0.00806$	&	$0.00802$	&	$0.01024$	\\
$0.55$	&	$-0.02702$	&	$-0.02575$	&	$-0.02590$	&	$-0.02537$	&	$1.60$	&	$0.01014$	&	$0.00642$	&	$0.00638$	&	$0.00553$	\\
$0.60$	&	$-0.02508$	&	$-0.02606$	&	$-0.02620$	&	$-0.02591$	&	$1.65$	&	$0.00797$	&	$0.01089$	&	$0.01086$	&	$0.00915$	\\
$0.65$	&	$-0.01250$	&	$-0.00779$	&	$-0.00781$	&	$-0.02359$	&	$1.70$	&	$0.00724$	&	$0.00896$	&	$0.00892$	&	$0.00826$	\\
$0.70$	&	$-0.01016$	&	$-0.00549$	&	$-0.00550$	&	$-0.00353$	&	$1.75$	&	$0.00828$	&	$0.00808$	&	$0.00805$	&	$0.00815$	\\
$0.75$	&	$-0.00412$	&	$-0.01408$	&	$-0.01417$	&	$-0.00428$	&	$1.80$	&	$0.01366$	&	$0.00658$	&	$0.00654$	&	$0.01361$	\\
$0.80$	&	$-0.00503$	&	$-0.00544$	&	$-0.00546$	&	$-0.00454$	&	$1.85$	&	$0.00776$	&	$0.00977$	&	$0.00974$	&	$0.01083$	\\
$0.85$	&	$-0.00505$	&	$-0.00656$	&	$-0.00658$	&	$-0.00505$	&	$1.90$	&	$0.01068$	&	$0.01346$	&	$0.01341$	&	$0.01192$	\\
$0.90$	&	$-0.00372$	&	$-0.00471$	&	$-0.00475$	&	$-0.00350$	&	$1.95$	&	$0.01537$	&	$0.01313$	&	$0.01305$	&	$0.01072$	\\
$0.95$	&	$-0.00503$	&	$-0.00452$	&	$-0.00454$	&	$-0.00406$	&	$2.00$	&	$0.01248$	&	$0.01225$	&	$0.01219$	&	$0.01031$	\\
$1.00$	&	$-0.00353$	&	$0.00024$	&	$0.00023$	&	$-0.00314$	&		&		&		&		&		\\

    \toprule[1.0pt]
    \end{tabular}
    \caption{Numerical results of the approximated LLE of SMIB system (\ref{eq:smib-c1})
    corresponding to the study-case 1, computed via the four $QR$-based techniques.}
    \label{tab:smib_c1}
\end{table}

\begin{table}[htb!]
    \centering
    \footnotesize
    \begin{tabular}{rrrrrrrrrr}
    \toprule[1.0pt]
    \multicolumn{1}{c}{$\rho$} & \multicolumn{1}{c}{LLE} &
    \multicolumn{1}{c}{$\rho$} & \multicolumn{1}{c}{LLE} &
    \multicolumn{1}{c}{$\rho$} & \multicolumn{1}{c}{LLE} &
    \multicolumn{1}{c}{$\rho$} & \multicolumn{1}{c}{LLE} &
    \multicolumn{1}{c}{$\rho$} & \multicolumn{1}{c}{LLE} \\
    \toprule[1.0pt]
    $0.00$	&	$-0.74586$	&	$0.60$	&	$-0.76208$	&	$1.20$	&	$-0.88585$	&	$1.80$	&	$-0.64744$	&	$2.40$	&	$-0.17751$	\\
$0.05$	&	$-0.74593$	&	$0.65$	&	$-0.76663$	&	$1.25$	&	$-0.89951$	&	$1.85$	&	$-0.64093$	&	$2.45$	&	$-0.09607$	\\
$0.10$	&	$-0.74572$	&	$0.70$	&	$-0.77893$	&	$1.30$	&	$-0.88447$	&	$1.90$	&	$-0.60895$	&	$2.50$	&	$-0.07736$	\\
$0.15$	&	$-0.74632$	&	$0.75$	&	$-0.78894$	&	$1.35$	&	$-0.89645$	&	$1.95$	&	$-0.56286$	&	$2.55$	&	$0.01194$	\\
$0.20$	&	$-0.74647$	&	$0.80$	&	$-0.78972$	&	$1.40$	&	$-0.87056$	&	$2.00$	&	$-0.53564$	&	$2.60$	&	$-0.00056$	\\
$0.25$	&	$-0.74694$	&	$0.85$	&	$-0.81366$	&	$1.45$	&	$-0.86469$	&	$2.05$	&	$-0.51171$	&	$2.65$	&	$0.04541$	\\
$0.30$	&	$-0.74786$	&	$0.90$	&	$-0.83571$	&	$1.50$	&	$-0.84152$	&	$2.10$	&	$-0.45690$	&	$2.70$	&	$0.16176$	\\
$0.35$	&	$-0.74901$	&	$0.95$	&	$-0.85522$	&	$1.55$	&	$-0.81798$	&	$2.15$	&	$-0.41292$	&	$2.75$	&	$0.20227$	\\
$0.40$	&	$-0.75015$	&	$1.00$	&	$-0.85872$	&	$1.60$	&	$-0.80238$	&	$2.20$	&	$-0.28496$	&	$2.80$	&	$0.27489$	\\
$0.45$	&	$-0.75187$	&	$1.05$	&	$-0.86683$	&	$1.65$	&	$-0.77636$	&	$2.25$	&	$-0.35850$	&	$2.85$	&	$0.28201$	\\
$0.50$	&	$-0.75370$	&	$1.10$	&	$-0.88458$	&	$1.70$	&	$-0.72155$	&	$2.30$	&	$-0.25838$	&	$2.90$	&	$0.33127$	\\
$0.55$	&	$-0.75980$	&	$1.15$	&	$-0.89293$	&	$1.75$	&	$-0.73500$	&	$2.35$	&	$-0.13307$	&	$2.95$	&	$0.38053$	\\

    \toprule[1.0pt]
    \end{tabular}
    \caption{Numerical results of the approximated LLE of SMIB system (\ref{eq:smib-c2})
    corresponding to the study-case 2, computed via \textit{C-EM} method.}
    \label{tab:smib_c2}
\end{table}

%%%%%%%%%%%%%%%%%%%%%%%%%%%%%%%% REFERENCES %%%%%%%%%%%%%%%%%%%%%%%%%%%%%%%%%%%%%%%%%%%%%%%%%%%%%%%
%\section*{Bibliography}
\bibliographystyle{abbrv}
\bibliography{references}
%%%%%%%%%%%%%%%%%%%%%%%%%%%%%%%%%%%%%%%%%%%%%%%%%%%%%%%%%%%%%%%%%%%%%%%%%%%%%%%%%%%%%%%%%%%%%%%%%%%
\end{document}